\documentstyle[12pt]{article}
\newtheorem{defin}{Definition}[section]
\newtheorem{th}[defin]{Theorem}
\newtheorem{prop}[defin]{Proposition}
\newtheorem{lem}[defin]{Lemma}
\newtheorem{cor}[defin]{Corollary}

\newcommand{\la}{\lambda }

\newcommand{\al}{\alpha }

\newcommand{\laal}{\mbox{$<\lambda,\alpha>$}}

\newcommand{\omeg}{\mbox{$\mid \! W\! \mid $}}

\newcommand{\ci}{\mbox{\bf{C}} }
\newcommand{\rr}{{\bf{R} }}
\newcommand{\pr}{\mbox{\bf{P}} } 
 \newcommand{\scp}{\mbox{ $<$ , $>$ } }
\newcommand{\cit}{\mbox{$\widetilde{C} $} } \newcommand{\ze}{{\bf{Z} }}
\newcommand{\pla}{\mbox{$P_{\lambda}$}}

\newcommand{\pro}{\mbox{{\em Proof.} } }
 \newcommand{\rim}{\mbox{{\em Remark} } }

\newcommand{\sal}{s_{\alpha}}
\newcommand{\nal}{n_{\alpha}}
\newcommand{\nalb}{n_{\alpha ,\be}}

\newcommand{\be}{\beta }
\newcommand{\ga}{\gamma }
\newcommand{\mi}{\mu_{i}}
\newcommand{\mj}{\mu_{j}}
\newcommand{\mh}{\mu_{h}}
\newcommand{\mii}{m_{i}}
\newcommand{\mij}{m_{j}}
\newcommand{\mih}{m_{h}}

\newcommand{\ej}{u_{j}}
\newcommand{\mth}{\tilde{\mu}_{h}}
\newcommand{\mti}{\tilde{\mu}_{i}}
\newcommand{\mtj}{\tilde{\mu}_{j}}
\newcommand{\deta}{(\eta)}
\newcommand{\dep}{(p)}
\newcommand{\seta}{(s_{\alpha}\eta)}

\newcommand{\thi}{t_{hi}}

\newcommand{\gtji}{\tilde{g}_{ji}}
\newcommand{\gthi}{\tilde{g}_{hi}}

\newcommand{\tih}{t_{h}}
\newcommand{\td}{\tilde}
\newcommand{\tj}{t_{j}}

\newcommand{\xa}{X_{\alpha}}
\newcommand{\xb}{X_{\beta}}
\newcommand{\xg}{X_{\gamma}}
\newcommand{\xma}{X_{-\alpha}}
\newcommand{\gi}{\mbox{\bf{g}} }
\newcommand{\tl}{\mbox{\bf{t}} }
\newcommand{\bl}{\mbox{\bf{b}} }
\newcommand{\ct}{{\cal T}}
\newcommand{\gl}{Gl(2)}
\newcommand{\pg}{PGl(2)}

\begin{document}
\begin{center} \mbox{ } \\ \mbox{ } \\ \mbox{ } \\ {\large
An elementary approach to the abelianization of the Hitchin system for 
arbitrary
reductive groups.} \\ \mbox{ } \\ by \\ Renata
Scognamillo%\footnote{partially supported by C.N.R.  and I.N.D.A.M.} \\
\mbox{
}
\\ {\em{Scuola Normale Superiore \\ Piazza dei Cavalieri, 7 \\ 56126 Pisa -
Italy} } \\ \mbox{ } \\ {\em{e.mail} }:  scognamillo@sabsns.sns.it \\
{\em{telefax} } 050-563513 \\ \mbox{ } \\ \end{center}
\section*{Introduction}
We consider here the  moduli space $\cal M$ of stable principal
$G$-bundles over a compact Riemann surface $C$, with $G$ an algebraic complex
group. We denote by $K$ the canonical bundle over $C$.  In   \cite{hi}
% em Stable bundles and integrable systems}, Duke Math.  J.
% bf{54} }(1987),
N.Hitchin  defined an analytic map $\cal H$ from the
cotangent bundle $T^{*}\cal M$ to the "characteristic space" $\cal K$
% with
% $dim\; {\cal K}=dim{\cal M}$
by associating to each $G$-bundle $P$ and section
$s\in H^{0}(C,adP\otimes K)$ the spectral invariants of $s$.
Hitchin showed for
$G=Gl(n),SO(n),Sp(n)$  that the generic fibre of $\cal H$ is an open set in an
abelian variety ${\cal A}$.  In fact, he considers in each case a
non-singular spectral curve $S$ covering $C$:  when $G=Gl(n)$, ${\cal A}$ is
identified with the Jacobian $J(S)$ ; in the other cases, there is a naturally
defined involution on $S$ and $\cal A$ is the associated Prym variety.  
%More
%recently, G.Faltings extended this theory to describe the moduli space of
%Higgs
%bundles for general reductive groups and proved
% by use of roots and weights
%that the generic fibre of the
%Hitchin map is a principal homogeneous space with respect to an
%algebraic group
% (namely the first \'{e}tale cohomology group of C with coefficients in a
% suitably defined group scheme)
%whose connected component of the identity is an abelian variety (see \cite{fa}
%).
More recently, Faltings extended these results and  described an  abelianization
 procedure for the moduli
space of Higgs $G$-bundles, with $G$
any reductive group (see \cite{fa} ). 
If $T\!  \subset \!
G$ is a fixed maximal torus with Weyl group $W$, one may construct for each
given generic element $\phi\in{\cal K}$ a ramified
covering \cit of $C$ having $\omeg $ sheets.  The combined action of $W$ on
\cit and on the group of one parameter subgroups of $T$ induces an action on
the space
of all principal $T$-bundles $\tau $ over \cit and we may
consider the subvariety $\widehat{\cal P}$
of those $\tau $ which are $W$-invariant in this sense. The connected
component ${\cal P}_{0}$ of $\widehat{\cal P}$ which contains
the trivial $T$-bundle is an
abelian variety.
In \cite{fa} it is shown that the generic
fibre of the Hitchin map is     a principal homogeneous space with respect to
a group
(namely the first \'{e}tale cohomology group of $C$ with coefficients in a
 suitably defined group scheme)
% (the first ....
% suitably defined in terms of \'{e}tale
% cohomology with
% coefficients in a group scheme,
which is isogenous to
$\widehat{\cal P}$.
% Here we approach the description of the Hitchin system
% ${\cal H}:T^{*}{\cal M}\rightarrow {\cal K}$ as follows.
%The aim of this paper is to describe explicitly the generic fibre of the
%Hitchin map
% ${\cal H}:T^{*}{\cal M}\rightarrow {\cal K}$ 
%in terms of a 
%"generalized Prym variety" $\cal P$.
In the present paper, by means of mostly elementary techniques, we 
%treat the  abelianization of the Hitchin system
% as follows.
 explicitly
construct a map ${\cal F}$ from each connected component
${\cal H}^{-1}(\phi )_{c}$ of ${\cal H}^{-1}(\phi )$
to ${\cal P}_{0}$ and show that $\cal F$ has finite fibres.  
%Our approach is the following.
 We use the classical theory of representations of
finite groups to compute $dim\; {\cal P}_{0}=dim\; {\cal M}$
and  
 conclude that the image under ${\cal F}$ of
% one connected component of the Hitchin fibre
${\cal H}^{-1}(\phi )$ contains a Zariski open set in ${\cal P}_{0}$. 

 In
case $G=PGl(2)$ one can check directly that the generic fibre of
${\cal F}_{c}: {\cal H}^{-1}(\phi )_{c}\rightarrow {\cal P}_{0}$ is a principal
homogeneous space with respect to a product of $(2\cdot deg\; K-2)$ copies of
$\ze
/2\ze$. However in case the Dynkin diagram of $G$ does not contain components 
of type $B_{l}$, $l\geq 1$ or 
when the commutator subgroup $(G,G)$ is simply connected  
% and for $G=SO(2n)$ 
the map ${\cal F}_{c}$ is injective.

Such results were announced in our previous paper \cite{sc} , in which we showed
that ${\cal P}_{0}$ is isogenous to a "spectral" Prym-Tjurin variety $\pla$ for
each given dominant weight $\la$.  
% \\
Results concerning the description of the Hitchin fibre in terms of 
generalized Prym 
varieties were also announced in R.Donagi, {\em Spectral covers}, 
preprint, alg-geom/9505009 (1995).

\section{The Hitchin map for any reductive group}
\label{sec-pre}
We denote by $C$ a compact Riemann surface of genus g$\; \geq 2$ and by $G$ a
reductive algebraic group over the field of complex numbers. We also write
$\gi $ as the Lie algebra of $G$.
The moduli space of stable principal $G$-bundles over $C$ is a
quasi-projective
 complex variety
$\cal M$ with $\; dim{\cal M}$=(g$-1)dimG+dimZ(G)$ , $Z(G)$
being the center
of $G$.
Note here that semistability for a principal $G$-bundle $P$ corresponds to
semistability for the holomorphic vector bundle $adP$ associated to the adjoint
representation \mbox{$Ad:G\rightarrow $gl$(\gi )\; $}
(\cite{ab} , \cite{ra} ).

We denote by $K$ the canonical line bundle over $C$.
By deformation theory and
Serre duality, a point in the cotangent bundle  $T^{*}\cal M$ of
$\cal M$ is a
pair $(P,s)$ with $P$ a stable principal $G$-bundle over $C$ and $s$ a section
of the vector bundle $adP\otimes K$.
The ring of  polynomials on $\gi $ which are invariant with respect to the
adjoint action is freely generated by homogeneous polynomials 
$h_{1},\ldots ,h_{k}$.
% Set  $deg\; h_{i}=d_{i}\; ,\; i=1,\ldots ,k$. 
Each $h_{i}$ induces a map ${\cal H}_{i}:adP\otimes K\rightarrow
K^{d_{i}}$ where $d_{i}= deg\; h_{i}$, 
and
the Hitchin map \[
{\cal H}\; :\; T^{*}{\cal M}\longrightarrow
{\cal K}={\displaystyle {\oplus_{i=1}^{k}}}H^{0}(C,K^{d_{i}})\]
takes $(P,s)$ to the element in $\cal K$ whose $i$-th component is the
composition of ${\cal H}_{i}$ with $s$ ( \cite{hi} ).
It is a remarkable fact that the dimension of $\cal K$
is equal to
the dimension of $\cal M$.
Moreover the map $\cal H$ is surjective. This fact can be deduced from the 
existence of very stable $G$-bundles (see \cite{lau}, \cite{br} ,\cite{kp} 
Lemma 1.4). 

%In \cite{hi} it is
% shown  for $G=Gl(n),Sp(n),SO(n)$ that the generic fibre of
%$\cal H$ is an open set in a certain abelian variety associated to a
% suitable
%"spectral curve" $S$ covering $C$.
We fix once and for all a maximal
torus $T\! \subset \! G$ with
associated root system $R=R(G,T)$ and Weyl group \mbox{$W\! =\! N_{G}(T)/T$}.
We also fix a subset $R^{+}\! \subset \! R$ (or equivalently a Borel subgroup
$B\!
\supset \! T$). If $\tl$ denotes the
Lie algebra of $T$, the differential of each root $\alpha \in R$ induces a map
\( d\alpha :\tl\otimes K\rightarrow K\) and the homogeneous
$W$-invariant
polynomials  
$\sigma_{1},\ldots ,\sigma_{k}$ on $\tl$ obtained by restriction of
$h_{1},\ldots ,h_{k}$ define a  Galois
covering \[ \underline{\sigma} =(\sigma_{1},\ldots ,\sigma_{k}):\tl
\otimes 
K\longrightarrow {\textstyle
\oplus_{i=1}^{k}}K^{d_{i}}\]
whose discriminant $\Xi$ is given by the zeroes of the
$W$-invariant function
$\prod _{\alpha\in R}d\alpha$ .
 For generic $\phi \in
{\cal K}\! =H^{0}(C,{\textstyle \oplus_{i}}K^{d_{i}})$, we consider the curve
$\widetilde{C}
:=\phi ^{*}(\tl\otimes K)$. This is a ramified covering of $C$ having
$\; m =\mid \! W\! \mid $ sheets, whose branch locus $Ram$
 satisfies by construction
\begin{equation} \label{eq:ram}
{\cal O}(Ram)\cong K^{\mid R\mid }\equiv K^{(dimG-rankG)}\; .\end{equation}
If we indicate by $\iota: \cit \rightarrow \tl\otimes K$ the natural inclusion
map, we have by definition, for each $w\in W$,
\begin{equation}
\label{eq:ad}
\iota(w\; \eta)\! =\! Ad(n_{w})\; \iota (\eta)
\end{equation}
where $n_{w}\! \in N_{G}(T)$  is any representative of $w$. Note also that, if
$\pi:\cit
\rightarrow C$ denotes the projection map, $d\alpha \circ \iota $ is a
holomorphic section of $\pi^{*}K$.
  \[ \begin{array}{rccc}
 & \widetilde{C}  & \stackrel{\iota}{\longrightarrow } & \tl\otimes K \\
   {\scriptstyle \pi} & \downarrow  &  &  \downarrow  \\
  &  C & \stackrel{\phi}{\longrightarrow } & {\textstyle \oplus _{i}}K^{d_{i}}
\end{array} \]
As a consequence of our genericity hypothesis, \cit has the following
properties: \\
 a) it is smooth and irreducible. \\
b)
each ramification point $p\in \pi^{-1}(Ram)$ has index 1; i.e. is a simple
zero
for the section
$\prod_{\alpha\in R^{+}}(d\alpha \circ \iota ):\cit \rightarrow \pi^{*}K^{\mid
R\mid /2}$.   \\[.2cm]
% Hence if $\widetilde{g} $ denotes the genus of
% $\widetilde{C} $ and $r$ the number of roots, we have
% \begin{eqnarray}
% 2-2\widetilde{g} & = & \omeg (2-2g)-r \frac{\omeg }{2} (2g-2)= \nonumber \\
%                 & = & \omeg (2-2g)(1+\frac{r}{2} )\; ;\label{eq:car} \\
%   \widetilde{g} & = & 1+\omeg (g-1)(1+\frac{r}{2} )\; . \nonumber
% \end{eqnarray}
This may be checked as follows. Let us denote by 
$\pi_{i}:K^{d_{i}}\rightarrow C$, $i=1,\ldots ,k$ and 
$q:\tl\otimes 
K\rightarrow C$ the projections.
%\begin{picture}
%\end{picture}
Moreover for every $i=1,\ldots ,k$  let us denote by
$\ga_{i}:K^{d_{i}}\rightarrow\pi_{i}^{*}K^{d_{i}}$ the tautological section.
For each $i$ we consider those sections of 
$q^{*}K^{d_{i}}$ that have the form 
$s=c\cdot\sigma_{i}^{*}\ga_{i}+q^{*}a_{i}$ for some $c\in\ci$ and $a_{i}\in  
H^{0}(C,K^{d_{i}})$. As $c$ varies in $\ci$ and $a_{i}$ in 
$H^{0}(C,K^{d_{i}})$ the zero divisor of $s$ forms a 
linear system 
$\delta_{i}$ of divisors in $\tl\otimes K$
that has no base 
points since the linear system $\mid K^{d_{i}}\mid $ on $C$
has no base points. 
For $\phi =(a_{1},\ldots ,a_{k})\in
{\cal K}$, the curve $\cit $ is defined by the equations
$\sigma_{i}^{*}\ga_{i}=q^{*}a_{i}$, $i=1,\ldots , k$. 
One immediately checks that the map 
\[
\begin{array}{lcc}
K^{d_{i}} & \longrightarrow & \pr ^{dim\; H^{0}(C,K^{d_{i}})}\\
     x   & \longmapsto  &  
[\ga_{i}(x),\pi_{i}^{*}a_{i,1}(x),\ldots ,\pi_{i}^{*}a_{i,m_{i}}(x)]
\end{array} \]
where the $a_{i,j}$'s form a basis of 
$H^{0}(C,K^{d_{i}})$
has image of dimension 2 and that
 $\sigma_{1}:\tl\otimes K\rightarrow K^{d_{1}}$ 
is dominant. By Bertini's theorem 
(see \cite{jou}, theorem 6.3) the 
divisor $X_{1}\in\delta_{1}$ of the section 
$\sigma_{1}^{*}(\ga_{1}-\pi_{1}^{*}a_{1})=\sigma_{1}^{*}\ga_{1}-q^{*}a_{1}$
with $a_{i}$ generic in
 $H^{0}(C,K^{d_{i}})$  
is smooth and irreducible. 
If $k\geq 2$, we next consider the linear system on $X_{1}$ given by the 
restriction of $\delta_{2}$. Since the polynomial  
$\sigma_{2}$ is 
algebraically independent from $\sigma_{1}$ the map
$\sigma_{2}\mid_{X_{1}}:X_{1}\rightarrow K^{d_{2}}$ is dominant. We use
the same argument as above and from
Bertini's theorem we obtain
that  the divisor $X_{2}\subset X_{1}$ of the section
$\sigma_{2}^{*}\ga_{2}-q^{*}a_{2}\mid_{X_{1}}$ with generic $a_{2}$
is smooth and irreducible. We can repeat the same argument for
the linear system 
$\delta_{i}\mid_{X_{i-1}}$ 
for every $i\leq k$ (since the map 
$\sigma_{i}\mid_{X_{i-1}}:X_{i-1}\rightarrow K^{d_{i}}$ is dominant)
and thus prove a) .
As for the statement b) one may consider the restriction of the linear 
systems above both to the discriminant locus $\Xi$ and to the locus
${\cal Z}\subset\Xi$ where  
$\prod _{\alpha\in R^{+}}d\alpha$ vanishes with multiplicity $\geq 2$ 
(${\cal Z}=Sing\; \Xi$). Again from Bertini's theorem one obtains that 
$\cit $ does not
 intersect $\cal Z$ and intersects $\Xi\setminus {\cal Z}$ transversely.
\\[.2cm] 
\rim 1.1.
For each $\al \in R^{+}$,
let $s_{\al }\in W$ denote the corresponding reflection.
As a consequence of condition b) above we may consider the
ramification locus in \cit as
a disjoint union: ${\cal D}=\coprod_{\alpha\in R^{+}}{\cal D}_{\alpha}$, with
${\cal D}_{\alpha}=\{ \mbox{zeroes of } \; \;  d\alpha \circ\iota \}
=\{ \eta\in \cit \mid  s_{\alpha}\; \eta\; =\eta\}$ . By our previous
considerations
${\cal D}_{\al }$ belongs to the linear system $\mid \! \pi^{*}K\! \mid $.
In case $G$ is simple and simply laced,
i.e. $W$ acts transitively on the set
of roots $R$,
we may write for each $y\in Ram$
\[
\pi^{-1}(y)=\! {\displaystyle
{\coprod _{\alpha \in R^{+}}}}{\cal D}_{\alpha }^{y}
\label{eq:div}
\]
where $\; {\cal D}_{\alpha }^{y}:={\cal D}_{\alpha }\cap
\pi^{-1}(y)$ is nonempty for every $\al \in R^{+}$.

If $G$ is not simply laced and has connected Dynkin diagram,
$R$ is the union of two $W$-orbits
$R_{1},\; R_{2}$,
each one consisting  of roots having the same length.
Then we have  \\
\begin{eqnarray}
\label{eq:orb}
\pi^{-1}(y)=\! \! {\displaystyle {\! \! \coprod _{\alpha \in
R_{1}\cap R^{+}}}}\! {\cal D}_{\alpha }^{y}
& \mbox{\makebox[1.5cm]{or} }
&
\pi^{-1}(y)=\! \! {\displaystyle {\! \! \coprod _{\alpha \in
R_{2}\cap R^{+}}}}\! {\cal D}_{\alpha }^{y}
 \end{eqnarray}
depending on whether  $y$ corresponds to a short or a long root.

More generally,
if the Dynkin diagram of $G$ has more than one connected component,
  we  have
as many different "kinds" of fibers
 \[ \pi^{-1}(y)=\! \! {\displaystyle {\! \! \coprod _{\alpha
\in  R_{j}\cap R^{+}}}}\! {\cal D}_{\alpha }^{y}\]
as are the $W$-orbits $R_{j}\! \subset \! R$.
Since for each $\al \in R^{+}$ we have $\mid \! {\cal D}_{\al }\! \mid =\omeg
\cdot deg\; K$ and each fibre over a branch point consists of $\omeg /2$
points, the number of fibres which correspond to the same orbit $R_{j}$ is
equal to
\begin{eqnarray}
n_{j}& = & \mid R_{j}^{+}\mid \cdot \omeg \cdot deg\; K/\frac{1}{2}
\omeg  \nonumber \\
\label{eq:nj}
 & = & \mid R_{j}\mid \cdot \; deg\; K\; .
\end{eqnarray}
\vspace{.3cm}

Let now $X(T)$  be the group of characters
on $T$ and
 consider the group $H^{1}(\widetilde{C},T)$ of
isomorphism classes of holomorphic principal
$T$-bundles  over $\widetilde{C}$.
 Each pair $(\tau ,\mu )$ with $\tau $ a principal $T$-bundle,
$\mu \in X(T)$
defines a line bundle $\tau _{\mu }\equiv
\tau \times _{\mu }\! {\bf C}$ and this way
$H^{1}(\widetilde{C},T)$ is identified with  \[
Pic(\widetilde{C})\otimes X(T)^{*},\]
 $X(T)^{*}\! \equiv \! Hom(X(T),\ze )$ being the dual group.
For the same reason,  the group of isomorphism classes of
topologically trivial principal
$T$-bundles   is a tensor product \[
J(\widetilde{C})\otimes X(T)^{*}\]
(here, as usual, $J(\cit )$ denotes the group
of divisors with  zero degree modulo  linear equivalence ).
Now, the action of $W$ on the sheets of $\widetilde{C}$
induces an action on $J(\cit )$. On the other hand, $W$ acts by
conjugation on $T$, hence on $X(T)^{*}$.
 If $\tau \! =D_{1}\otimes \chi_{1}+\cdots
+D_{l}\otimes \chi_{l}$ is a principal
$T$-bundle over $\widetilde{C}$ and $w\in W$ an element of the Weyl group,
we set \[
^{w}\tau =w\; D_{1}\otimes \; ^{w}\chi _{1}+\cdots
+w\; D_{l}\otimes \; ^{w}\chi_{l}\; .\]
\begin{defin}
The generalized Prym variety
$\; {\cal P}=[J(\widetilde{C})\otimes X(T)^{*}]^{W}$
consists of those isomorphism classes of topologically trivial $T$-bundles
$\tau $ which satisfy $\; ^{w}\tau \equiv \tau \; $ for each
$w\in W$.
\end{defin}
Note that
$\cal P$ is an algebraic group whose connected component of
the identity ${\cal P}_{0}$ is an abelian variety.
\section{Computing the dimension of $\cal P$}
The following can be deduced from the above mentioned
Faltings' result describing the generic Hitchin
fibre as isogenous to 
$\widehat{\cal P}=[Pic(\widetilde{C})\otimes X(T)^{*}]^{W}$ 
(\cite{fa}, theorem III.2)
and the fact (due to G.Laumon and proved in \cite{fa}, theorem II.5) that 
all Hitchin fibers have the same dimension:   
\begin{prop} \label{it:prd}
The dimension of $\cal P$ is equal to the dimension of $\cal M$.
\end{prop}
% \pro
In this section we want to give a direct proof of such statement.
If we set ${\cal S}\equiv \! X(T)\otimes_{\ze }\ci $ and denote by $H^{1}$ the
first
cohomology $W$-representation $H^{1}(\cit ,\ci )$, by Doulbault theorem we have
\[ dim{\cal P}=\frac{1}{2} dim[H^{1}\otimes {\cal S}^{*}]^{W}=\frac{1}{2}
dim\; Hom_{W}({\cal S},H^{1}).\]
We will compute
 \( M\equiv dim\; Hom_{W}({\cal S},H^{1})\)
% =<\chi_{{\cal S}},\chi_{H^{1}}>\)
by use of the
classical theory of representations of finite groups and associated characters
(for more details about this subject, see for example \cite{se} ).
% We recall here the main points of such theory.

For any $W$-representation $V$ considered here, 
% with $F$ some finite group, 
we
 denote by
$\chi_{V}:W\rightarrow \ci $ its character 
(for $\rho :W\rightarrow Gl(V)$ the
homomorphism defining the representation, we have by definition
$\chi_{V}(w)=trace(\rho(w))\; \forall w\in W$).
By the theory of characters of finite groups
we have 
\begin{equation} \label{eq:car}
 M=<\chi_{{\cal S}},\chi_{H^{1}}>\end{equation}
where $\scp $ is the usual scalar product between characters.
If $N$ is the number of connected components of the Dynkin diagram $\Pi$ of
$G$ and $h=dimZ(G)$ we have a decomposition \[
{\cal S}=\underbrace{{\cal B}\oplus \cdots \oplus {\cal B}}_{h} \oplus {\cal
S}_{1}\oplus \cdots \oplus {\cal S}_{N}\]
where $\cal B$ is the 1-dimensional trivial representation and ${\cal S}_{i}$
the irreducible reflection representation corresponding to the $i$-th component
of $\Pi $, $i=1,...,N$. 
%By property {\bf 1} above we have 
Then we may rewrite (\ref{eq:car}) as
\begin{equation}
\label{eq:dim}
M=h<\chi_{{\cal B}},\chi_{H^{1}}>+\sum_{i=1}^{N}<\chi_{{\cal
S}_{i}},\chi_{H^{1}}>. \end{equation}
We observe that $W$ acts trivially on the cohomology groups
$H^{0}(\cit ,\ci )\cong H^{2}(\cit ,\ci )\cong \ci $.
Hence the
Lefschetz character $\chi_{L}\equiv \chi_{H^{0}}-\chi_{H^{1}}+\chi_{H^{2}}$
satisfies $\chi_{L}=2\chi_{{\cal B}}-\chi_{H^{1}}$ and we have
\begin{eqnarray}
\label{eq:lef}
<\chi_{{\cal B}},\chi_{H^{1}}>  & = & 2-<\chi_{{\cal B}},\chi_{L}> \\
\label{eq:lefbis}
<\chi_{{\cal S}_{i}},\chi_{H^{1}}> & = &
-<\chi_{{\cal S}_{i}},\chi_{L}>.\end{eqnarray}
On the other hand, it is well known (Hopf trace formula, see e.g.\cite{cr} )
that the Lefschetz character satisfies \[
\chi_{L}=\chi_{\widetilde{C}^{0}}-\chi_{\widetilde{C}^{1}}+\chi
_{\widetilde{C}^{2}}\]
$\cit ^{n}$ being the free $\ci $-module generated by the $n$-cells of some
cellular decomposition of $\cit $
($\cit ^{n}\cong H_{n}(K^{n},K^{n-1};\ci )$, with $K^{j}$ the $j$-th
skeleton of $\cit $, $j=n,n-1$).

We choose one finite triangulation $\Delta$ of $C$ whose set of
vertices
contains all branch points. We denote by  $C^{n}$ the free module generated
by the $n$-cells  of $\Delta$ for $n=1,2$ , and by $C_{0}^{0}$ and
$D_{j}$ the free modules whose generators are respectively
all vertices not lying
in the branch locus $Ram$ and all branch points corresponding to the same
$W$-orbit $R_{j}\! \subset \! R$ (see {\em Remark} 1.1.).
Let $N'$ be the number of $W$-orbits in $R$, and for each $j=1,\ldots N'$
let us fix one positive root
$\al _{j}\! \in \! R_{j}^{+}$  and set $H_{j}=\{
1,s_{\al _{j}}\} \! \subset \! W$. We denote by
$Ind^{W}_{H_{j}}(B_{j})$ the $W$-representation
induced by the 1-dimensional trivial representation $B_{j}$ of
$H_{j}$
 (by definition ,$Ind^{W}_{H_{j}}(B_{j})=\oplus_{[w]\in W/H_{j}}\ci
 v_{[w]}$ with $W$ acting by $u\circ v_{[w]}=v_{[uw]}$).
We have the
following isomorphisms of $W$-modules:
\begin{eqnarray*} \cit ^{2} & \cong  & \ci [W]\otimes C^{2} \\
\cit ^{1} & \cong   & \ci [W]\otimes C^{1} \\
\cit ^{0} & \cong   & \ci [W]\otimes C^{0}_{0}\oplus
\bigoplus_{j=1}^{N'}Ind^{W}_{H_{j}}(B_{j})\otimes
D_{j} \\
  &  \equiv & \ci [W]\otimes C^{0}_{0}\oplus
\bigoplus_{j=1}^{N'}(Ind^{W}_{H_{j}}(B_{j}))^{n_{j}}\end{eqnarray*}
where $\ci [W]$ denotes as usual the regular representation and the
$n_{j}$'s are defined as in (\ref{eq:nj}) . \\
By Frobenius reciprocity formula 
we have \[
<\chi_{\cal B},\chi_{Ind^{W}_{H_{j}}(B_{j})}>=<\chi_{B_{j}},\chi_{B_{j}}>=1\; 
;\]
% hence by property {\bf 4} 
and since from the general theory
each irreducible $W$-representation occurs as a subrepresentation
of $\ci [W]$ as many times as is its dimension, 
we obtain
\begin{equation} \label{eq:ada}
 <\chi_{{\cal B}},\chi_{L}>=rk\; C^{2}-rk\; C^{1}+rk\; C^{0}_{0}+\mid
Ram\mid =(2-2g).\end{equation}
% Again by properties {\bf 4}  and {\bf 5} we get
% Using again such property for $\ci [W]$ and Frobenius reciprocity formula 
Analogously, we have
\[ <\chi_{{\cal
S}_{i}},\chi_{L}>=(rk\; C^{2}-rk\; C^{1}+rk\;
C^{0}_{0})dim{\cal S}_{i}+\sum_{j=1}^{N'}n_{j}<\chi_{B_{j}},\chi
_{res_{j}{\cal S}_{i}}>\]
where 
$res_{j}{\cal S}_{i}$ denotes 
the representation obtained by restriction
% S}_{i}$  
to  $H_{j}$.\\
Now, given some positive root $\alpha\in R^{+}$, the corresponding
reflection $s_{\alpha}\in W$ acts trivially on ${\cal S}_{i}$ whenever
$\alpha\! \notin \! {\cal S}_{i}$ , otherwise
% , in case $\alpha\in {\cal R}_{i}$,
 it acts trivially on one subspace of codimension 1 in
${\cal S}_{i}$. Thus
we get
\begin{eqnarray}
<\chi_{{\cal S}_{i}},\chi_{L}> & = & (rk\; C^{2}-rk\; C^{1}+rk\;
C^{0}_{0})dim{\cal S}_{i}+\sum_{R_{j}\subset {\cal S}_{i}}
n_{j}(dim{\cal S}_{i}-1)+\nonumber \\
 &  + &
\sum_{R_{j}\not\subset {\cal S}_{i}}
n_{j}\cdot dim{\cal S}_{i}\nonumber \\
\label{eq:adb}
 & = &  (2-2g)\; dim{\cal S}_{i}-\sum_{R_{j}\subset {\cal S}_{i}}
n_{j}\; .
\end{eqnarray}
By substituting (\ref{eq:ada})  and (\ref{eq:adb}) respectively in
(\ref{eq:lef}) and (\ref{eq:lefbis})  and then (\ref{eq:lef})
and (\ref{eq:lefbis})  in (\ref{eq:dim}) , we finally obtain
\begin{eqnarray*}
 M & = & 2h+(2g-2)(h+\sum_{i=1}^{N}dim{\cal S}_{i})+\sum_{j=1}^{N'}
n_{j} \\  & = & 2h+(2g-2)dim\; T+\mid
\! Ram\! \mid  .\end{eqnarray*}
Since $dim\; T+\mid \! R\! \mid =dim\; G$,
by (\ref{eq:ram}) we get \[
dim{\cal P}\equiv \frac{1}{2} M=(g-1)dim\; G+h.
% \makebox[3cm][r]{$\Box$}
\]
\section{The main results}
\label{sec-res}
In this section we will define a map $\cal F$ from each component of the
generic Hitchin fibre to the abelian variety ${\cal P}_{0}$ and study its
properties. We first show how one can associate to each given pair
$(P,s)\in {\cal H}^{-1}(\phi)$ a $T$-bundle ${\cal T}={\cal T}(P,s)$ which
satisfies $^{w}\! \ct \cong \ct \; \forall w\in W$.

For $\phi\in {\cal K}$ generic, let then $P$ be a principal $G$-bundle and
$s\in H^{0}(C,adP\otimes K)$ such that $(P,s)\in {\cal H}^{-1}(\phi)$.
We first consider the restriction $P_{0}$ of $P$ to the open set $C_{0}$.
Since for every $\xi\in C_{0}$, $s(\xi)\in \gi$ is regular semisimple 
(for an analysis of the regular elements in $\gi$ , 
%with respect to the adjoint representation
 see for example  \cite{ko}),
we have a morphism of vector bundles
\[ [s,\hspace{.2cm} ]:adP_{0}\longrightarrow adP_{0}\otimes K\]
whose kernel $\cal N$ is a bundle of Cartan subalgebras in $\gi$. 
We thus have
a section \[
\gamma :C_{0}\rightarrow P/N_{G}(T)\equiv P\times_{G} G/N_{G}(T)\]
 locally defined by
$\gamma(\xi)=\nu(\xi)N_{G}(T)$ where $\nu(\xi)\in G$ satisfies 
$Ad\; \nu(\xi)\tl ={\cal N}_{\xi}\equiv c_{\gi}(s(\xi))$. 
If we pull back $P_{0}$ over $\cit_{0}$
we actually have  a section
\begin{equation}\label{eq:sect}
 \varphi :\cit_{0}\rightarrow \pi^{*}P_{0}/T
\end{equation}
locally defined by
$\varphi\deta=\mu\deta T$ where  $\mu\deta\in G$ satisfies 
\begin{equation}\label{eq:diag}
Ad\;\mu\deta (\iota\deta)=s(\pi\deta).\end{equation}
Thus over $\cit_{0}$ the bundle $\pi^{*}P$ has a reduction of its 
structure group to
$T$. Moreover, from (\ref{eq:ad}) we have for each $w\in W$
\begin{equation}
\label{eq:nw}
\varphi(w\eta)=\mu\deta n_{w}^{-1}T\end{equation}
which implies that such $T$-reduction
$\tau_{0}=\varphi^{*}(\pi^{*}P_{0})$ 
is $W$-invariant with respect to the action
previously defined. 
Now if we consider
  a Borel subgroup $B\subset G$ containing
$T$, 
%which corresponds to our choice of $R^{+}$ 
the inclusion map $T\hookrightarrow B$ and $\varphi$
define a section 
$:\cit_{0}\rightarrow\pi^{*}P\times_{G} G/B$.
Since $G/B$ is a complete variety, by the valuative criterion of properness 
this section
can be extended to the whole curve $\cit $ and we thus obtain
 (uniquely up to isomorphisms) 
a  $B$-reduction $P_{B}$ of the $G$-bundle $\pi^{*}P$ such that
$P_{B}\mid_{\scriptsize \cit_{0}} $ is the $B$-extension of 
% $T$-bundle 
 $\tau_{0}$.\\[.1cm]
If $(\; ,\; )$ denotes a $W$-invariant scalar product on $X(T)_{\ze}\otimes 
\rr$ and $\be\in R$,
we define as usual the one parameter subgroup
$\be'\in Hom(X(T),\ze)$ by \begin{equation} \label{eq:'}
\be'(\lambda)=<\la,\be >\equiv
\frac{2(\la ,\be )}{(\be ,\be )}\; \;
\forall \lambda\in X(T).\end{equation}
We want to prove the following:
\begin{th} \label{it:thw}
Let $\tau_{B} =\tau(P,s)$ be the $T$-bundle
associated to 
%the $B$- bundle 
$P_{B}$
via the natural projection
$B\rightarrow T$. Let us fix one theta characteristic $\frac{1}{2} K$ and
consider the $T$-bundle
$K_{\rho}=\frac{1}{2}\pi^{*}K\otimes \sum_{\be\in R^{+}} \be'$, where
$R^{+}\subset R$ is the subset of positive roots that corresponds to $B$.
Then  ${\cal T}(P,s):=\tau_{B} +K_{\rho}$ is $W$-invariant.
\end{th}
The proof will be organized in a few lemmas.  
We
first observe that since $W$ is generated by the simple
reflections it suffices to show \begin{equation}
\label{eq:tbun}
 ^{s_{\al}}\tau_{B} \cong \tau_{B} +\pi^{*}K\otimes \al '
\end{equation}
for every simple root $\al $. In fact
we have ${\displaystyle \sum_{\be \in R^{+}}
s_{\al }(\be ')=\sum
_{\stackrel{\be\in R^{+}}{\be \neq \alpha}} \be ' -\al '} $,
so, if
relation (\ref{eq:tbun}) holds, one has
$^{s_{\al}}(\tau_{B} +K_{\rho})\cong \tau_{B} +K_{\rho}$.
In terms of line bundles associated to characters on $T$, relation
(\ref{eq:tbun}) can be rewritten as
\begin{equation}
\label{eq:fon}
%s_{\al}^{*}(\tau_{B}\times_{\sal(\la)}\ci)\cong \tau_{B}\times_{\la}\ci
(^{s_{\al}}\tau_{B} -\tau_{B})\times_{\la}\ci
\cong \laal\pi^{*}K\; 
\; \;
\forall \la \in X(T).
\end{equation}
Given a simple root $\al$, let us denote by $\sal(B)$ the Borel subgroup
$n_{\al}Bn_{\al}^{-1}$, where $\nal\in N_{G}(T)$ represents $\sal$.
One analogously obtains another $T$-bundle $\tau_{\sal (B)}$ such that
$\tau_{\sal (B)}\mid_{\scriptsize\widetilde{C}_{0}}\cong\tau_{0}$ from the 
completion of
$\tau_{0}$ to 
%the section $\varphi$ in (\ref{eq:sect}) to  
%$\hat{\varphi}_{\sal(B)}:\cit\rightarrow \pi^{*}P/\sal(B)$ and the 
%corresponding 
an $\sal(B)$-reduction $P_{\sal(B)}$.  
The first lemma treats the relationship between
 $\tau_{B}$ and 
$\tau_{\sal (B)}$.
\begin{lem}\label{it:borels}
We have $\tau_{\sal(B)}\cong\; ^{\sal}\tau_{B}$.
\end{lem}
\pro
We consider an open covering $\{ V_{h}\} _{h\in H}$ of $C$
over which $P$ and the canonical bundle $K$ can be
trivialized
and  with the property that each $V_{h}$
contains at most one branch point.  We choose a \v{C}ech covering
${\cal U}=\{ U_{h}\}_{h\in H}$
 of \cit  to be given by  all open sets $U_{h}=\pi^{-1}(V_{h})$
(by definition each $U_{h}$ is stable with respect to the action of $W$).
For $h\in H$ we choose frames $e_{1}^{h},\ldots ,e_{q}^{h}$
 for the vector
bundle $adP\otimes K$ over $V_{h}\subset
C$, $q$ being equal to the dimension
of $\gi $. With respect to this
choice the section
 $s:C\rightarrow ad P\otimes K$ is locally given
by "coordinates" $s_{h}:V_{h}\rightarrow \gi$
satisfying
\begin{equation} \label{eq:shl}
s_{h}=Ad\; g_{hl}\cdot
k_{hl}s_{l}\; \; \mbox{ for }
V_{h}\cap V_{l}\neq \emptyset ,
\end{equation}
$g_{hl}$ and $k_{hl}$ being transition functions for $P$, $K$
respectively. Let $\iota_{h}:U_{h}\rightarrow \tl $ be coordinates for
$\iota:\cit\rightarrow\tl\otimes K$.
We define $J\subset H$ to be the subset
of those indices $j$ such that $V_{j}$ contains a branch point and set
 $I=H\setminus J$.
For each $h\in H$ we fix maps $\mh:U_{h}\rightarrow G$ such that
for each $i\in I$  $\mi$ satisfies 
\begin{equation}\label{eq:mui}
Ad\;\mi\deta (\iota_{i}\deta)=s_{i}(\pi\deta)\end{equation}
(compare with (\ref{eq:diag}) ) 
and the 0-chain $\{ \mh\deta B\}_{h\in H}$ defines the section 
$\widehat{\varphi}_{B}:\cit\rightarrow\pi^{*}P/B$ completing $\varphi$
in (\ref{eq:sect}).
By definition, the $B$-bundle $P_{B}$ is represented by the cocycle
$\{ b_{hl}\} \in{\cal Z}^{1}({\cal U},B)$ where
$b_{hl}\deta \equiv
\mh \deta ^{-1}g_{hl}(\pi\deta )\mu_{l} \deta$. 
Define
$\{ b'_{hl}\} \in{\cal Z}^{1}({\cal U},\sal(B))$
by $b'_{hl}\deta=\nal b_{hl}(\sal\eta)n_{\al}^{-1}\; \forall\eta\in 
U_{h}\cap U_{l}$.
We have
$b'_{hl}\deta \equiv
\nal\mh \seta ^{-1}g_{hl}(\pi\deta )\mu_{l} \seta\nal^{-1}$, hence
$\{ b'_{hl}\}$ represents an $\sal (B)$-reduction of $\pi^{*}P$.
On the other hand, from (\ref{eq:nw}) we have 
$\{ \mu_{i} \seta\nal^{-1}T\}_{i\in I}=\{ \mu_{i} \deta T\}_{i\in I}$
hence $\{ b'_{hl}\} $ represents $P_{\sal(B)}$.
Now, if we denote by $p:B\rightarrow T,\; 
p':\sal(B)\rightarrow T$ the natural projections we have
$p'\circ b'_{hl}\deta =\nal (p\circ b_{hl}(\sal\eta))n_{\al}^{-1}$
(since every Borel subgroup is a semidirect product of its maximal torus and 
its maximal unipotent subgroup). 
Since $\{ \nal (p\circ b_{hl}(\sal\eta))n_{\al}^{-1}\} $ are by definition
transition functions for $^{\sal}\tau_{B}$, we thus have an isomorphism
$\tau_{\sal(B)}\cong\; ^{\sal}\!\tau_{B}$.\hspace{.2cm} $\Box$\\
%Let
%${\cal U}=\{ U_{h}\}_{h\in H}$
% be as above. 
We keep the notations of the proof of lemma \ref{it:borels}.
For each positive root $\beta\in R^{+}$, we
shall denote by $\beta_{h}:U_{h}\rightarrow \ci$ the coordinates of the
section of $\pi^{*}K$ over \cit given by the composition $d\beta \circ \iota$
(see \S \ref{sec-pre}). 
 Our next step consists in finding suitable transition
functions $b_{ji}$ for $P_{B}$ on intersections $U_{i}\cap U_{j}$
with $j\in J$.
Indeed, we will find suitable maps 
$\mu_{j}:U_{j}\rightarrow G$ with $j\in J$ defining the completed section 
$\widehat{\varphi}_{B}$. 
We fix  nilpotent generators 
$\{ \xg\}_{\ga\in R^{+}}$ in the Lie algebra $\bl $ 
of $B$
with $ad\; t(\xg )=\gamma (t)\xg$ $\forall t\in \tl $, $\; \forall\ga\in 
R^{+}$.  
In general, the completion $\hat{\varphi}_{B}:\cit\rightarrow \pi^{*}P/B$ 
of our $\varphi$ above %(see (\ref{eq:sect}))
is locally given by  
holomorphic maps $f_{j}:U_{j}\rightarrow G$ with $j\in J$
such that
\begin{equation}
\label{eq:ff}
Ad\; f_{j} \deta ^{-1}s_{j}(\pi\deta)=\iota_{j} \deta +\sum_{\gamma\in R^{+}}
a_{\gamma}\deta \xg\; .
\end{equation}
By \rim 1.1 , for $j\in J$ the set $U_{j}$ is a union
of open sets $\bigcup_{\beta\in R(j)\cap R^{+}}U_{j,\beta}$
%  with
% $U_{j,\beta}\cap \pi^{-1}(Ram)\subset {\cal D_{\beta}}$
where $R(j)$ is some $W$-orbit of  roots  depending on $j$ and
each $U_{j,\beta}$ contains only those
ramification points that are zeroes for $\beta_{j}$ . 
\begin{lem} \label{it:bor}
There exists a holomorphic map $\mj :U_{j}\rightarrow G$
satisfying for each $\be\in R(j)\cap R^{+}$ and $\eta\in U_{j,\be}$
\begin{equation}
\label{eq:jb}
Ad\; \mj \deta ^{-1}s_{j}(\pi\deta)=\iota_{j} \deta +\xb \; . 
\end{equation}
%where $\xb^{(r)} \in \bl$ is constant
%on each connected component $U_{j,\be}^{(r)}\subset U_{j,\beta}$.
\end{lem}
\pro
We construct $\mj$ separately on each open set  $U_{j,\be}$. 
By our genericity hypothesis we may assume for
a ramification point $p\in U_{j,\be}$
\begin{equation}
\label{eq:fff}
Ad\; f_{j} \dep ^{-1}s_{j}(\pi\dep )=\iota_{j} \dep + \xb
\end{equation}
with 
% $a_{\beta}\dep \neq 0$ and 
$\beta_{j}\dep \equiv
d\beta(\iota_{j}\dep )=0$.\\
Let $\alpha$ be the root with minimal height in
$R^{+}\setminus\{ \beta\} $
such that $a_{\alpha}\deta $ in (\ref{eq:ff}) is not identically zero.
The map
$c_{j}\deta
  =exp\frac{a_{\alpha}\deta }{\alpha_{j}\deta }\xa :U_{j,\be}\rightarrow G$ is
holomorphic on each fixed connected component of $U_{j,\beta}$ and
by evaluating $Ad\; c_{j}\deta $ on the right-hand side of (\ref{eq:ff})
we get
 \[ Ad\; c_{j}\deta (\iota_{j} \deta +
\sum_{\gamma\in R^{+}}
 a_{\gamma}\deta \xg )
=\iota_{j} \deta + a'_{\beta}\deta \xb +
\sum
_{\stackrel{\gamma\in R^{+}\setminus \{ \be\} }{\gamma >\alpha}}
 a_{\gamma}\deta \xg \; .\]
By an induction argument we can then 
%rewrite (\ref{eq:ff} )
assume
\begin{equation}
% \label{eq:ffff}
Ad\; f_{j} \deta ^{-1}s_{j}(\pi\deta )=\iota_{j} \deta +a_{\beta}\deta \xb 
\end{equation}
% ($f_{j}$ having been suitably renamed).
where $a_{\be}(p)=1$ (since we may multiply $f_{j}$ by a suitable constant in 
$T$).
Consider now the map 
$d_{j}\deta
  =exp\frac{a_{\beta}\deta
-1}{\beta_{j}\deta }\xb :U_{j}\rightarrow G$. 
Since $p$ is a simple zero for $\beta_{j}\deta $, $d_{j}\deta $ is
holomorphic on each chosen connected component of $U_{j,\beta}$. 
We have  \[
Ad\; d_{j}\deta (\iota_{j} \deta +a_{\beta}\deta 
\xb )=\iota_{j} \deta + \xb \]
and the claim of our lemma is proved.
%(with $\xb^{(r)}= $).
\makebox[1.2cm][r]{$\Box$} \\
%As a consequence of this lemma, for $j\in J$ such that
%$U_{i}\cap U_{j}\neq \emptyset $ 
%we have for each $\be\in R(j)$ and $\eta\in 
%U_{j,\beta}\cap U_{i}$
%\begin{equation} \label{eq:e}
%Ad\; ( \mj \deta ^{-1}g_{ji}\deta \mi \deta )\;
%k_{ji}\deta \iota_{i}\deta
%=\iota_{j}\deta +\xb .\end{equation}
For each  $j\in J$,
define $\ej :U_{j}%\cap U_{i}
\rightarrow B$ by
 $u_{j}\deta =exp\;
\frac{\xb }{\be_{j}\deta }$  whenever $\eta\in U_{j,\be}$. We have 
\begin{equation} \label{eq:uj}
Ad\; \ej\deta^{-1}\iota_{j}\deta=\iota_{j}\deta+\xb\; .\end{equation}  
We may represent the completed section $\widehat{\varphi}_{B}$
by $\{ \mh\deta B\}$ where  the $\mi $'s are as in (\ref{eq:mui})
for every $i\in I$ and the $\mj $'s satisfy  
(\ref{eq:jb}) for every $j\in J$. By substituting  
(\ref{eq:mui}) and (\ref{eq:jb}) in (\ref{eq:shl})
and replacing $\iota_{j}\deta+\xb $ with $Ad\; \ej\deta^{-1}\iota_{j}\deta $
we obtain transition functions on each nonempty intersection $U_{j}\cap U_{i}$ 
\begin{equation}
\label{eq:gp}
 b_{ji}\deta \equiv
\mj \deta ^{-1}g_{ji}(\pi\deta )\mi \deta =u_{j}^{-1}\deta
t_{ji}\deta
\end{equation}
where  $t_{ji}\deta :U_{i}\cap U_{j}\rightarrow T$ is holomorphic (as 
$u_{j}$ is holomorphic on $U_{i}\cap U_{j}$ ).
Since each element in $B$ can be written uniquely as a product of a unipotent 
element by an element in $T$ we have
%$\mi=\mu\mid_{U_{i}}$ and
$t_{ji}=p\circ b_{ji}$.
%:U_{i}\cap U_{j}\rightarrow T$.
%We remark that choosing another set
%$\{ X'_{\ga}\}_{\ga\in R^{+}}$ 
%of elements in $\bl$
%with 
%$ad\; t\; X'_{\ga}=\gamma (t)X'_{\ga}$ $\forall t\in \tl $ 
%$\forall\ga\in 
%R^{+}$ 
%  corresponds to multiplying 
%each $t_{ji}\mid_{U_{i}\cap U_{j,\ga}}$ by a suitable element in 
%$T$.

We now compare $P_{B}$ with $P_{\sal(B)}$. By definition we only need to
compare them around the ramification points. 
As set of nilpotent generators
in the Lie algebra of $\sal(B)$ we may choose 
$\{ \xb\} _{\be\in R^{+}\setminus\{\al\} }\cup\{ Ad\; \nal (\xa)\}$.
Thus from lemma \ref{it:bor}
we may define a section
$\hat{\varphi}_{\sal(B)}:\cit\rightarrow \pi^{*}P/\sal(B)$ completing 
$\varphi$ by 
\begin{eqnarray*}
\hat{\varphi}_{\sal(B)}\deta & = & \mj\deta\sal(B)\;
\makebox{ for }\; \eta\in U_{j}\setminus U_{j,\al}  \\
\hat{\varphi}_{\sal(B)}\deta & = & \mj\seta\nal^{-1}\sal(B)\;
\makebox{ for }\; \eta\in U_{j,\al} 
\end{eqnarray*} 
where the $G$-valued maps $\mj$ satisfy 
(\ref{eq:jb}). From this we see that
 $P_{\sal(B)}$ and $P_{B}$ are isomorphic 
on $\cit\setminus {\cal D}_{\al}$ and that on all intersection sets 
$U_{j,\al}\cap U_{i}$ with $j\in J$ we have transition functions for
$P_{\sal(B)}$ of the form 
\begin{equation}\label{eq:b'}
b'_{ji}\deta =\nal\mj\seta ^{-1}\mj\deta b_{ji}\deta .
\end{equation}
% for each $\eta\in U_{j,\al}$, $\; j\in J$.\\
%On the other hand 
If we apply lemma \ref{it:bor} to the set $\sal(R^{+})$ of positive roots
corresponding to $\sal(B)$ 
%(note
%that the decompositions ${\cal D}=\bigcup_{\be\in \sal(R^{+})}{\cal D}_{\be}$ 
%and $U_{j}=\bigcup_{\be\in \sal(R^{+})\cap R(j)}U_{j,\be}$ remain 
%the same)
we obtain on $U_{j,\al}\cap U_{i}$  a factorization 
$b'_{ji}\deta ={u'_{j}}^{-1}\deta t'_{ji}\deta $ 
with $u'_{j}\deta =exp\;
\frac{Ad\; \nal(\xa ) }{-\al_{j}\deta }=\nal {u_{j}}^{-1}\deta\nal^{-1}$ and 
$t'_{ji}\deta =p'\circ b'_{ji}\deta $ (compare with (\ref{eq:gp}) ).
Let us denote by I the identity element in $G$.
From (\ref{eq:b'}) and lemma \ref{it:borels}  
a meromorphic 
section of $\; ^{\sal}\tau_{B} -\tau_{B}$ is  given by a 0-cochain
$\{\tih \}_{h\in H}\in{\cal C}^{0}({\cal U},T)$ where
\begin{eqnarray}
 \tih\deta & = & \mbox{I~} \mbox{ whenever } h\in I \mbox{ or } h\in J 
\mbox{ and } 
\eta\notin U_{j,\al}
%outside } \bigcup_{i\in I}U_{i}\cup\bigcup_{j\in J}U_{j,\al} 
\label{eq:sezi}\\
\tj\deta & = & \nal u_{j}\deta^{-1}\mj\seta ^{-1}\mj\deta u_{j}\deta^{-1} 
\; \; \forall\eta\in U_{j,\al},\; j\in J.   \label{eq:sez}
\end{eqnarray}
By (\ref{eq:jb}) on each $U_{j,\al }$ the map   
$h_{j}\deta =\mj\seta ^{-1}\mj\deta$ 
satisfies
\begin{equation}\label{eq:Fj}
Ad\; h_{j}\deta
(\iota_{j}\deta +\xa )=\iota_{j}\seta +\xa 
=Ad\; \nal(\iota_{j}\deta )+\xa  .  
\end{equation}
%(where we have set
Choose $\xma \in \gi$ so that $\xa ,\xma ,h_{\al
}:=[\xa 
,\xma ]\in \tl $ generate a Lie subalgebra ${\bf h}_{\al }\subset
\! \gi $
with ${\bf h}_{\al }\cong sl(2)$ and
$d\al (h_{\al })=2$. Define
\[ F_{j}\deta =exp(\al _{j}\deta \xma )\; \; \forall \eta\in
U_{j,\al }.\] Since $F_{j}\deta$ satisfies 
$Ad\; F_{j}\deta
(\iota_{j}\deta +\xa )=Ad\; \nal(\iota_{j}\deta )+\xa $,
  by (\ref{eq:Fj}) we have
on $U_{j,\al}$
\begin{equation}
\label{eq:c}
 \mj \seta ^{-1}\mj \deta =F_{j}\deta \cdot L_{j}\deta
\end{equation} 
%or equivalently
%\begin{equation} \label{eq:eq9}
%\mj \seta ^{-1} =F_{j,\al}\deta \cdot L_{j,\al}\deta\cdot \mj \deta^{-1}
%\end{equation}
where for each  $\eta\in U_{j,\al}$, $L_{j}\deta \in B$ 
lies in the centralizer of $\iota_{j}\deta
+\xa \in \bl $.  Note that for  $q$
any ramification point in $U_{j,\al }$ we have by definition
\begin{equation} \label{eq:uno}
L_{j}(q)=\mbox{I} .\end{equation}
In particular the map $L_{j}$ is
holomorphic. 
%We define  for each $\eta\in U_{j,\al}$
%By the definition of $L_{j}$ and
%(\ref{eq:uno}) one can easily check that 
Since when  $\eta\in U_{j,\al}$
%\setminus {\cal D}_{\al}$  
is not a ramification point
$~\iota_{j}\deta + \xa$
is regular semisimple and by (\ref{eq:uj}) one has 
$c_{\gi}(\iota_{j}\deta + \xa)=Ad\; 
u_{j}\deta^{-1}\tl $ ,
the holomorphic $T$-valued 
map $l_{j}\deta =p\circ L_{j}\deta$ has the form 
\begin{equation}
\label{eq:lj}
l_{j}\deta =u_{j}\deta L_{j}\deta
u_{j}\deta
^{-1}.\end{equation}
Relation (\ref{eq:sez}) becomes
\begin{equation}\label{eq:tjz}
\tj\deta =z_{j}\deta\cdot l_{j}\deta
\end{equation}
where the map 
$z_{j}\deta\equiv \nal u_{j}\deta^{-1} F_{j}\deta u_{j}\deta ^{-1}$
has values in $T$ and is holomorphic
everywhere in $U_{j,\al}$ but on the ramification points.
 The connected subgroup $H_{\al }\subset \! G$ generated by
$exp(\xa ),exp(\xma ),exp(h_{\al })$  is isomorphic
to a copy of $Sl(2)$ or $PGl(2)$ in $G$ and
one can compute $z_{j}\deta $ 
%by (\ref{eq:zj}) 
directly in terms of two by two
matrices. In the $Sl(2)$ case, denoting by $\varrho$ the isomorphism
$:H_{\al}\rightarrow Sl(2)$, one has for some $c\in\ci^{*}$
\begin{eqnarray}
\nonumber
\varrho (z_{j}\deta ) & = &
\mp\left( \begin{array}{cc}
0 & -1 \\ 
1 & 0
\end{array}\right)
\left( \begin{array}{cc}
1 & -c/\al_{j}\deta \\ 
0 & 1
\end{array}\right)
\left( \begin{array}{cc}
1 & 0 \\ 
\al_{j}\deta /c & 1
\end{array}\right)
\left( \begin{array}{cc}
1 & -c/\al_{j}\deta \\ 
0 & 1
\end{array}\right) \\
\label{eq:cru}
 & = &
k\cdot
diag(\al _{j}\deta ,\al _{j}\deta ^{-1})
%\left( \begin{array}{cc}
%\al_{j}\deta & 0 \\ 
%0 & \al_{j}\deta ^{-1}
%\end{array}\right)
\end{eqnarray}
where $k\in T$ is a constant and $\al _{j}\deta $
%is defined as %in \ref{eq:bj}.
are the coordinates of the section $d\al\circ\iota$, according to our previous
notations.  As for $H_{\al}\stackrel{\varrho}{\cong}\pg$ one gets
\begin{equation}
\label{eq:crup}
% \al (z_{j}\deta )=\al _{j}\deta ^{-2}$.\\
\varrho (z_{j}\deta )=
\overline{k\cdot diag(\al _{j}\deta ,\al _{j}\deta ^{-1})}
\end{equation}
where 
the bar indicates the image under the factor map $:\gl\rightarrow\pg$.
Let now $T_{\al }\subset \! T$ be the identity component of the subgroup
$Ker(\al)=\{ t\in T\; \mid \; \al (t)=1\} $. The centralizer
$Z_{\al}$ in $G$ of $T_{\al}$ is a reductive group of semisimple rank 1 having
Lie algebra ${\bf z}=\tl \oplus \ci \xa \oplus \ci \xma $, and
it is known that such a group is a product $T'\times H$ , $T'$ being
a torus and H being a copy of $Sl(2)$, $PGl(2)$ or $Gl(2)$.
The case $H=Sl(2)$ is characterized by the group of characters $X(T)$ being
an orthogonal direct sum $\ze \chi_{1} \oplus X'$, with $\chi_{1}=\sqrt{\al}$.
If we compose 
any 
 $\la \in X'$ 
% (we have $<\la ,\al >=0$), composing $\la $ 
with the 0-chain
$\{ t_{h}\} _{h\in H}$
% in
% (\ref{eq:d}), (\ref{eq:dd}), (\ref{eq:ddd})
defined by (\ref{eq:sezi}) and (\ref{eq:sez})
we obtain a nowhere vanishing holomorphic section of
the line bundle
$(^{\sal }\tau_{B}-\tau_{B})\times _{\la}\ci $.
If instead we compose $\chi_{1}$ to  
$\{ t_{h}\} _{h\in H}$, by 
 (\ref{eq:tjz}) and (\ref{eq:cru}) 
% (we have $<\la ,\al >=1$),
we get an holomorphic section
for
$(^{\sal }\tau_{B}-\tau_{B})\times _{\chi_{1}}\ci $
having  simple zeroes exactly on the locus
  ${\cal D_{\al}}$.
%=\{ p\in \pi^{-1}(Ram)\; \mid\; d\al \circ \iota
%\dep =0\}$. 
Thus relation (\ref{eq:fon}) is satisfied (see \rim 1.1).

The case $H=PGl(2)$ is characterized by $X(T)$ being
an orthogonal direct sum $\ze \al \oplus X'$.
For $\la \in X'$, we get the same result as for the $Sl(2)$ case.
% composing $\la $ to the functions
% $\{ t_{h}\} _{h\in H}$ gives us a never vanishing holomorphic section of
% the line bundle
% $\sal ^{*}(\tau_{\la}^{-1})\otimes \tau_{\la}$.\\
For $\la =\al $ we find instead an holomorphic section
for
$(^{\sal }\tau_{B}-\tau_{B})\times _{\la}\ci $
having  zeroes of
multiplicity
two on ${\cal D_{\al}}$
% =\{ p\in \pi^{-1}Ram\; \mid\; \al \circ \iota
% \dep =0\}$;
. This proves (\ref{eq:fon}).

In case $H=Gl(2)$, we have an orthogonal direct sum
$X(T)=X'\oplus \ze \chi_{1}\oplus \ze \chi_{2}$ with
$\al =\chi_{1}\cdot \chi_{2}^{-1}$.
Composing $\la \in X'$ gives us again
$^{\sal }\tau_{B}\times _{\la}\ci \cong \tau_{B}\times _{\la}\ci $
 as in the previous cases.
 If we compose $\chi_{1}$ we
obtain an holomorphic section of
$(^{\sal }\tau_{B}-\tau_{B})\times _{\chi_{1}}\ci $
having simple zeroes exactly on
${\cal D_{\al}}$.
If we compose $\chi_{2}$ we
obtain a meromorphic section of
$(^{\sal }\tau_{B}-\tau_{B})\times _{\chi_{2}}\ci $
having simple poles exactly on
${\cal D_{\al}}$.
Thus relation (\ref{eq:fon}) holds also in this case and theorem \ref{it:thw}
is proved.
%\makebox[2cm][r]{$\Box$} \\[.2cm]

We thus have a map \[ \begin{array}{cccc}
  {\cal T}: & {\cal H}^{-1}(\phi) &
 \rightarrow  &  \widehat{{\cal P}} \equiv [Pic(\cit 
)\otimes X(T)^{*}]^{W} \nonumber \\
 & (P,s)  & \longmapsto & \tau (P,s)+K_{\rho}\; \; .
\end{array} \]
Note that from (\ref{eq:tbun}) and lemma \ref{it:borels} 
\hspace{.1cm} $\cal T$ does not 
depend on the choice of the Borel subgroup $B\supset T$ (or of the subset
of positive roots in $R$). 
\begin{defin}
% \label{it:defF}
Let ${\cal H}^{-1}(\phi)_{c}$
%\subset {\cal H}^{-1}(\phi)$
be some
 connected component of ${\cal H}^{-1}(\phi)$. For
a fixed point $(P',s')\in {\cal
H}^{-1}(\phi)_{c}$ we define ${\cal F}_{c}:{\cal
H}^{-1}(\phi)_{c}\rightarrow {\cal P}_{0}$
by
% \begin{equation} \label{eq:F}
\[
 {\cal
F}_{c}(P,s)={\cal T}(P,s)-{\cal T}(P',s')\equiv {\tau(P,s)}-{\tau(P',s')}\;
. % \subset {\cal P}_{0}
\]
% ${\cal T}$ and $\tau $ being defined as in theorem {\em \ref{it:thw}}.
\end{defin}
Such  definition does not depend on our previous choice of the
theta characteristic $\frac{1}{2} K$. We now
want to study the fibers of ${\cal
F}_{c}$. First we make the following \\
{\em Remark} \ref{sec-res}.1
\hspace{.2cm} For $i\in I$, the maps $\mi \deta$ in
%(\ref{eq:gp})
(\ref{eq:mui}) are defined up to multiplication to the right by some holomorphic
map $m_{i}:U_{i}\rightarrow T$. As for $j\in J$, any other holomorphic
map $\mj '\deta$ satisfying (\ref{eq:jb}) 
has the form $\mj '\deta =\mj
\deta
M_{j}\deta $ where for every $\al\in R(j)\cap R^{+}$,
$M_{j}\deta :U_{j,\al }\rightarrow B$ is holomorphic
and such that $M_{j}\deta \in c_{G}(\iota_{j}\deta +\xa )$. If we
replace $\mj$ and $\mi$ 
%in (\ref{eq:gp}) 
with
the new maps $\mj '\deta $ and $\mi '\deta =\mi \deta m_{i}\deta $, we
obtain
from $(P,s)$ and $B$ 
%$P_{B}$ 
an equivalent cocycle $\{ m_{h}^{-1}\thi m_{i}\} $ representing
$\tau_{B}$.
Since for every $j\in J$ and $q\in U_{j}\cap {\cal D}_{\al}$ 
$~\iota_{j}(q)+\xa \in \bl$ is regular, we have 
$c_{G}(\iota_{j}(q)+\xa )=T_{\al}{\cal U}_{\al}$, 
where $T_{\al}$ is the identity component of
$Ker(\al:T\rightarrow\ci^{*})$ and ${\cal U}_{\al}$ is the 
unipotent 1-dimensional subgroup corresponding to the root $\al $. 
Hence
% The
%fact that $M_{j}\deta $ is holomorphic implies that
%for each $j\in J$
the
$T$-valued map $m_{j}\deta :=p\circ M_{j}\deta
\equiv u_{j}\deta M_{j}\deta u_{j}\deta ^{-1}$
satisfies for every $\al\in R(j)\cap R^{+}$ 
 \begin{equation} \label{eq:norm}
\al (m_{j}(q))=1\;  \mbox{ $\forall q\in U_{j}\cap {\cal D}_{\al}$ .} 
 \end{equation}
\begin{lem} \label{it:l2}
Let $(P,s),(Q,v)$ be pairs in ${\cal
H}^{-1}(\phi)$ such that $\tau(P,s)$ and $\tau(Q,v)$ are isomorphic.
Let $\{ t_{hl}\}$ and $\{ \td{t}_{hl}\}$ with $h,l\in  H$ be cocycles 
representing $\tau(P,s)$ and $\tau(Q,v)$ respectively and suppose
\begin{equation}
\td{t}_{hl}=\mih^{-1}t_{hl}m_{l} \label{eq:t1p} 
%\mtj^{-1}\gtji\mti=\ej^{-1}\mij^{-1}\tji\mii \nonumber % \label{eq:gpp}
\end{equation}
where the maps $m_{h}:U_{h}\rightarrow T$ are holomorphic and satisfy
condition {\em (\ref{eq:norm})} for
every $j\in J$ and $\al\in R(j)\cap R^{+}$. 
%satisfy the condition $\al(\mij(p))=1$ whenever
%$p$ lies in ${\cal D}_{\al }\cap U_{j}$ .
Then $Q$ is isomorphic to $P$ and $v=s$.
\end{lem}
 \pro  
For what concerns  $P$ and the construction of $\tau (P,s)$ we keep the 
notations used in the proof of theorem \ref{it:thw} .
In particular we still consider a C\v{e}ch covering 
${\cal U}=\{ U_{h}\}_{h\in H}$ of $\cit $ consisting of $W$-invariant 
open sets as it was first defined in the proof of lemma \ref{it:borels}. 
For each nonempty intersection $U_{h}\cap U_{l}$ we have 
transition functions for
the $B$-reduction $Q_{B}$ of $\pi^{*}Q$ having the form:
\begin{eqnarray}\label{eq:gpt}
\td{b}_{ji}\deta 
& = & \mtj\deta^{-1}\gtji(\pi\deta )\mti\deta =\ej\deta^{-1}\td{t}_{ji}\deta
\; \; \forall j\in J,i\in I \\
\td{b}_{hi}\deta 
& = & \mth\deta^{-1}\gthi(\pi\deta )\mti\deta =\td{t}_{hi}\deta
\; \; \forall i,h\in I \end{eqnarray}
where $\{ \td{g}_{hl}\}_{h,l\in H}$ are transition functions for the $G$-bundle
$Q$ and $\mti$, $\mtj$ are defined 
analogously as $\mi$ and $\mj$ in (\ref{eq:gp}). 
%By substituting (\ref{eq:t1}) in (\ref{eq:t1p}) we get
%\begin{equation}
%\mth^{-1}\gthi\mti=\mih^{-1}\mh^{-1}\ghi\; \mi\; \mii\; \; \; \forall h,i\in I.
%\label{eq:l1}
%\end{equation}
For $j\in J$, define $M_{j}:U_{j}\rightarrow B$ by \begin{equation}
\label{eq:mj}
M_{j}:=\ej^{-1}\mij\ej \; \; \; 
\mbox{ (see \rim \ref{sec-res}.1 ).} \end{equation}
The hypothesis of the lemma provide that $M_{j}$ is holomorphic on  
 $U_{j,\al}$ for each $\al\in R(j)\cap R^{+}$ and  we have $M_{j}\deta
\in c_{G}(\iota_{j}\deta +\xa )\; \forall \eta \in U_{j,\al }$ by
definition of $u_{j}$.
%Thus, by (\ref{eq:gp}), we have
%\begin{eqnarray}
%\label{eq:l2}
%M_{j}^{-1}\mj^{-1}\gji\mi\mii & = & M_{j}^{-1}\ej^{-1}\tji\; \mii \nonumber \\
% & = & E_{j}^{-1}\mij^{-1}\tji\; \mii \nonumber \\
% & = & \mtj^{-1}\gtji\; \mti .
%\end{eqnarray}
Define the holomorphic maps 
% $\Gamma_{l}:U_{l}\rightarrow G$ by
\begin{eqnarray*}
\Gamma_{i} & = & \mi\; \mii\; \mti^{-1}\; \; \forall i\in I 
%\label{eq:Gi} 
\mbox{ and} \\
\Gamma_{j} & = & \mj M_{j}\mtj^{-1}\; \; \forall j\in J.
%\label{eq:Gj} 
\end{eqnarray*}
From (\ref{eq:gpt}), (\ref{eq:gp}) and (\ref{eq:t1p})
we obtain the equivalence condition
between cocycles on $\cit$:
\[
\tilde{g}_{hl}(\pi\deta )=\Gamma_{h}\deta ^{-1}g_{hl}(\pi\deta )\Gamma_{l}\deta
\; \; \; \forall\eta\in U_{h}\cap U_{l}\; \forall h,l\in H.\]
The claim of the lemma is then proved provided we show that the
maps $\Gamma_{l}$ are invariant with respect to the action of $W$
on the sheets of $\cit$.
In fact if we indicate by $\{ v_{h}\} _{h\in H}$ the coordinates of $v$
so that $v_{h}=Ad\tilde{g}_{hl}\cdot k_{hl}v_{l}$,
by our definition of the maps $\tilde{\mu}_{l}$, $\tilde{\mu}_{h}$
%(compare with relation (\ref{eq:jb}) ) 
we have:\[
Ad\; \Gamma_{l}\; v_{l}=s_{l}\; \; \forall l\in H.\]
%the $\{ Y_{l}\} _{l\in H}$ being coordinates for the section $s$.
Since $W$ is generated by the simple reflections, it suffices to show
$\Gamma_{l}\seta =\Gamma_{l}\deta$ for every simple reflection $\sal$. 
From (\ref{eq:nw})
we have for each $i\in I$
\begin{equation}\label{eq:mis}
{\mu}_{i}\seta ^{-1}{\mu}_{i} \deta =\nal {l}_{i}\deta \end{equation}
for suitable holomorphic maps $l_{i}:U_{i}\rightarrow T$. By evaluating
the transition functions
$t_{hi}=\mh^{-1}g_{hi}\mi$ with $h,i\in I$ on $\sal \eta$ and 
replacing $\mi\seta$ with 
${\mu}_{i} \deta {l}_{i}\deta^{-1}\nal^{-1}$ and
$\mh\seta$ with 
${\mu}_{h} \deta {l}_{h}\deta^{-1}\nal^{-1}$  we obtain  
\begin{equation} \label{eq:t}
{t}_{hi}\seta =\nal {l}_{h}\deta {t}_{hi}\deta
{l}_{i}\deta^{-1} \nal^{-1}.\end{equation}
Analogously, if we define $\td{l}_{i}:U_{i}\rightarrow T$ by 
\begin{equation}\label{eq:mtis}
\mti\seta^{-1}\mti\deta=\nal\td{l}_{i}\deta,
\end{equation} we have
\begin{equation} \label{eq:td}
\td{t}_{hi}\seta =\nal \td{l}_{h}\deta \td{t}_{hi}\deta
\td{l}_{i}\deta^{-1} \nal^{-1}.\end{equation}
By replacing $\td{t}_{hi}$ with $\mih^{-1}t_{hi}\mii$ in both sides of
(\ref{eq:td}) and substituting (\ref{eq:t}) in the left-hand side, we
obtain 
an equality both sides of which contain only factors with values in $T$.
We cancel  $t_{hi}\deta$ and obtain
\[  m_{h}\deta\cdot  
\nal^{-1} m_{h}\seta^{-1} \nal\cdot\td{l}_{h}\deta^{-1}\cdot l_{h}\deta
=m_{i}\deta\cdot  
\nal^{-1} m_{i}\seta^{-1} \nal\cdot\td{l}_{i}\deta^{-1}\cdot l_{i}\deta
\]
for every $\eta\in U_{h}\cap U_{i}$, $i,h\in I$.
We can repeat the same calculation on intersection sets
$U_{i}\cap U_{j}$ with $j\in J$ and $i\in I$. What we need is the analog
for $j\in J$ of the relations (\ref{eq:mis}) and (\ref{eq:mtis}).
On each open set $U_{j,\al}$ the map $\mj\deta$ is related 
with $\mj\seta$ via the identity (\ref{eq:c}).
If for each $\be\in R^{+}\setminus\{ \al\}$
we define $n_{\al\be}\in N(T)$ to be the representative of $\sal$ satisfying 
$Ad\; n_{\al,\be}(\xb )=X_{\sal(\be)}$, by construction of the maps 
$\mj$ in lemma (\ref{it:bor}) 
we have for $\eta\in U_{j,\be}$
\begin{equation}\label{eq:mjs}
{\mu}_{j}\seta ^{-1}{\mu}_{j} \deta = n_{\al,\be}{L}_{j}\deta 
\end{equation} 
where $L_{j}\deta$ is a suitable element 
in the centralizer of $\iota_{j}\deta +\xb $.
%and satisfies $L_{j}(q)=\mbox{I} $ for each $q\in {\cal D}_{\be}$ 
%(compare with (\ref{eq:uno}) ).
We analogously
define  
$\td{L}_{j}:U_{j}\rightarrow B$ $\forall j\in J$ by
%\( \left\{
\begin{eqnarray} \label{eq:ctd}
\td{\mu}_{j}\seta ^{-1}\td{\mu}_{j} \deta & = & F_{j}\deta \td{L}_{j}\deta 
\; \; \mbox{ for $\eta\in U_{j,\al}$} \\
\label{eq:mtjs}
\td{\mu}_{j}\seta ^{-1}\td{\mu}_{j} \deta & = & \nalb\td{L}_{j}\deta 
\; \; \mbox{ for $\eta\in U_{j,\be}$ with $\be\neq\al$ }
\end{eqnarray}
and set for each $\eta\in U_{j}$
\begin{eqnarray}
\label{eq:plj}
{l}_{j}\deta &:= & p\circ L_{j}\deta =u_{j}\deta {L}_{j}\deta
u_{j}\deta
^{-1}\\
\label{eq:pltj}
\td{l}_{j}\deta & := & p\circ \td{L}_{j}\deta =u_{j}\deta \td{L}_{j}\deta
u_{j}\deta
^{-1}.\end{eqnarray}
One uses (\ref{eq:c}), (\ref{eq:ctd}) and the fact that the map
$z_{j}\deta=\nal u_{j}^{-1}\deta F_{j}\deta u_{j}^{-1}\deta$ 
(see (\ref{eq:tjz}) ) is holomorphic $T$-valued
outside the ramification points (hence it commutes
with any 
other map with values in $T$), to obtain
by the same procedure described above for all pairs of indices $h,i\in I$
\[  m_{j}\deta\cdot  
\nal^{-1} m_{j}\seta^{-1} \nal\cdot\td{l}_{j}\deta^{-1}\cdot l_{j}\deta
=m_{i}\deta\cdot  
\nal^{-1} m_{i}\seta^{-1} \nal\cdot\td{l}_{i}\deta^{-1}\cdot l_{i}\deta
\]
for each $\eta\in U_{j,\al}\cap U_{i}$.
One uses (\ref{eq:mjs}) and (\ref{eq:mtjs}) to prove the same identity
for all $\eta\in U_{j,\be}\cap U_{i}$ with $\be\neq\al$.
In conclusion, the maps  \\
$m_{h}\deta\cdot 
\nal^{-1} m_{h}\seta^{-1} \nal\cdot
\td{l}_{h}\deta^{-1}\cdot l_{h}\deta:U_{h}\rightarrow T$  with  $h\in H$
are the
restriction to $U_{h}$ of a global holomorphic map on $\cit$, hence are
equal to some constant ${\bf c}$. We compute such map
on
one ramification point $q\in U_{j,\al }$ . Since we have 
 $l_{j}(q)=\td{l}_{j}(q)=\mbox{I} $ (compare with (\ref{eq:uno}) )
and  $\al (\mij (q))=1$ by hypothesis, we obtain
${\bf c}=$I ,
i.e. \begin{equation}
\label{eq:lm}
m_{h}\seta = \nal m_{h}\deta \cdot l_{h}\deta\cdot \td{l}_{h}\deta ^{-1}
\nal^{-1}\; \; \; \; 
\forall
h\in H. \end{equation}
By use of (\ref{eq:mis}), (\ref{eq:mtis}) and this last
identity we find $\Gamma_{i}\seta =\Gamma_{i}\deta$ for each $\eta\in U_{i}$,
$i\in I$.  
%(see (\ref{eq:Gi}) ), 
%we have \\[.2cm]
%$\Gamma_{i}\seta
%=\mi\seta\mii\seta \td{\mu}_{i}\seta^{-1}=$ \vspace{.2cm} \\
%               $=\mi\deta
%l_{i}\deta^{-1}\nal^{-1}\mii\seta\nal\td{l}_{i}\deta\td{\mu}_{i}\deta^{-1}
As for $j\in J$, if $\eta$ is in $U_{j,\al}$ we have 
by  (\ref{eq:c}) and (\ref{eq:ctd}), by the definition of $M_{j}$, $l_{j}$ and
$\td{l}_{j}$ and by (\ref{eq:lm}) 
%\\[.2cm]
\begin{eqnarray*}
\Gamma_{j}\seta & = &
\mi\deta u_{j}\deta^{-1}
l_{j}\deta^{-1}z_{j}\deta^{-1}\mij\deta l_{j}\deta\td{l}_{j}\deta^{-1}
z_{j}\deta\td{l}_{j}\deta u_{j}\deta\td{\mu}_{j}\deta^{-1}=\\
%[.1cm]
%\hspace{1cm} 
& = & \Gamma_{j}\deta .
\end{eqnarray*}
%\\[.2cm]
If $\eta$ is in $U_{j,\be}$, one proves \( \Gamma_{j}\seta =\Gamma_{j}\deta \)
by using (\ref{eq:mjs}), (\ref{eq:mtjs}), (\ref{eq:lm}) and the identity
(following from the above definition of $n_{\al,\be}$) 
 $n_{\al,\be}\; u_{j}\seta\; n_{\al,\be}^{-1}=u_{j}\deta$. 
\makebox[2cm][r]{$\Box$}
\begin{lem}
\label{it:l3}
Let $(P,s),(Q,v)$ be pairs in ${\cal
H}^{-1}(\phi)$ such that $\tau(P,s)$ and $\tau(Q,v)$ are isomorphic.
Let $\{ t_{hl}\}$ and $\{ \td{t}_{hl}\}$ with $h,l\in  H$ be cocycles 
representing $\tau(P,s)$ and $\tau(Q,v)$ respectively and write
\begin{equation}
\td{t}_{hl}=\mih^{-1}t_{hl}m_{l} %\label{eq:t1p} 
%\mtj^{-1}\gtji\mti=\ej^{-1}\mij^{-1}\tji\mii \nonumber % \label{eq:gpp}
\end{equation}
for suitable holomorphic maps
 $m_{h}:U_{h}\rightarrow T$ with $h\in H$.
Up to multiplying each $m_{h}$ by one and the same suitably chosen element 
in $T$,
the following holds:
\\ {\em (i)}  for each positive root $\al\in R^{+}$ and $q\in U_{j}\cap {\cal
D}_{\al}$ we have $\al(\mij(q))=\mp 1$.
{\em (ii)}   if for $\al\in R^{+}$ there exists some character $\la\in X(T)$
such that \begin{equation} \label{eq:sc}
<\la,\al>=1\; ,
\end{equation}
we have $\al(\mij(q))=1$ $\forall q\in U_{j}\cap {\cal
D}_{\al}$.
\end{lem}
\pro 
Choose one ramification point $q_{\al}\in {\cal D}_{\al}$ for each $\al\in 
\Delta $, $q_{\al}\in U_{j(\al )}$ for suitable $j(\al)\in J$. Up to 
multiplying the maps $\{ m_{h}\}_{h\in H}$ by a suitable element in 
$T$ we may assume
\begin{equation} \label{eq:ga}
\al(m_{j(\al)}(q_{\al}))=1\; \; \forall \al\in \Delta.
\end{equation} 
We keep the same notation as before.
We consider the maps $\{ l_{h}\} $ and $\{ \td{l}_{h}\} $, $h\in H$ as 
in (\ref{eq:mis}), (\ref{eq:mtis}), (\ref{eq:plj}) and (\ref{eq:pltj})  
% Let us denote by $\td{t}_{ij}$ the cocycle $\mij^{-1}\tji\mii :U_{j}\cap
% U_{i}\rightarrow T$ with $j\in J,i\in I$,
 and 
let $\al$ be some simple root.
From the proof of lemma (\ref{it:l2}) one has 
that the maps
$m_{h}\deta\cdot 
\nal^{-1} m_{h}\seta^{-1} \nal\cdot
\td{l}_{h}\deta^{-1}\cdot l_{h}\deta:U_{h}\rightarrow T$ 
are the restriction of a global holomorphic map on $\cit$. Computing 
such map on $q_{\al}$ gives us by 
(\ref{eq:ga}) and the fact that we have
%$l_{j(\al )}(q_{\al})=\td{l}_{j(\al )}(q_{\al})=\mbox{I} 
$l_{j}(q)=\td{l}_{j}(q)=\mbox{I} 
\; \; \forall q\in {\cal D}_{\al}\cap U_{j}
$ 
\begin{equation} \label{eq:s}
m_{j}(q)\cdot \nal^{-1}m_{j}(\sal q)^{-1}\nal
\cdot \td{l}_{j}(q)^{-1}\cdot l_{j}(q)=\mbox{I} \; \; \forall q\in 
{\cal D}\cap U_{j}, j\in J
\end{equation}
and
\[
% \label{eq:m}
 \mij(q)=\; \; \nal^{-1}m_{j}(\sal q)\nal\; \; \forall q\in 
{\cal D}_{\al}\cap U_{j},\; j\in J.
\]
 By evaluating $\al:T\rightarrow\ci^{*}$ on both sides of
this last identity we obtain \[
\al^{2}(\mij(q))=1.\]
If moreover $\al$ satisfies condition (\ref{eq:sc}), evaluating $\la$ on
both sides of the same identity 
%(\ref{eq:m}) 
gives \(
\la(\mij(q))=\la(\mij(q))\cdot\al^{-1}(m_{j}(q))\),
or \[ \al(\mij(q))=1.\]
The claim of the theorem is thus proved for every simple root.
Consider now $q\in {\cal D}_{\be}$ with $\be\in R^{+}\setminus\Delta$. 
Note that for $q\in U_{j}$, from the definition of $l_{j}$ and $\td{l}_{j}$
and the fact that $L_{j}(q)$ and  $\td{L}_{j}(q)$ belong to the centralizer
in $G$ of $\iota_{j}(q)+\xb $
we have
\begin{equation}\label{eq:bl}
\be(l_{j}(q))=\be(\td{l}_{j}(q))=1
\end{equation}
(compare with (\ref{eq:norm}) in {\em Remark} 3.1).
By evaluating  $\be :T\rightarrow\ci^{*}$ on both sides of
(\ref{eq:s}) as  $\al$ runs over all simple 
roots we obtain
$\be (\mij(q))=\be(\nal^{-1}m_{j}(\sal q)\nal )\; 
\; $ $\forall\al\in\Delta$, hence \[
\be(\mij(q))=\be(n_{w}^{-1}m_{j}(wq)n_{w})\; \; \forall w\in W.\]
On the other hand, we know that there exist $\al\in\Delta$ and $u\in W$ 
with $u(\al)=\be$. We thus have
%$\mij(u^{-1}(q))=n_{u}^{-1}m_{j}(q)n_{u}$ we obtain \[
\[
\be(\mij(q))=\be(n_{u}m_{j}(u^{-1}q)n_{u}^{-1})=
\al(\mij(u^{-1}q))=\mp 1.\mbox{\makebox[2cm][r]{$\Box$} } \]
% By substting relation (\ref{eq:dd}) to the right-hand side of the identity
% \[ \td{t}_{ji}\seta=\mij^{-1}\seta\tji\seta\mii\seta\]
% we get \\
% $\td{t}_{ji}\seta
% =\mij ^{-1}\seta
% z_{j}\deta({^{\sal}l_{j}}\deta)({^{\sal }t}_{ji}\deta)({^{\sal
% }l_{i}}\deta)^{-1}\mii\seta $= \vspace{.2cm} \\
%  $= \mij ^{-1}\seta
% z_{j}\deta
% ({^{\sal}l_{j}}\deta)({^{\sal}m_{j}}\deta)({^{\sal
% }\td{t}_{ji}}\deta)({^{\sal}m_{i}}\deta)^{-1}({^{\sal
% }l_{i}}\deta)^{-1}\mii\seta
% .$ \\
% But now  relation (\ref{eq:cru})
% \rim 2.1
%  holds also for the cocycle  $\td{t}_{ji}$ ;
% thus for any ramification point $p\in U_{j}\cap {\cal D}_{\al}$ we have \[
% \mij(p)^{-1}\cdot {^{\sal}l_{j}}(p)\cdot {^{\sal}m_{j}}(p)=I.\]
% On the other hand, by the same relation we know that ${^{\sal}l_{j}}(p)=I$,
% hence we find
\begin{th} \label{it:cor1}
Suppose $G$ has one of the following properties: \\
a) the commutator group $(G,G)$ is simply connected ;\\
b) the Dynkin diagram of $G$ has no component of type $B_{l},\; l\geq 1$.\\
 Then 
% for every connected
% component ${\cal H}^{-1}(\phi)_{c}\subset {\cal
% H}^{-1}(\phi)$ the map ${\cal F}_{c}:{\cal
% H}^{-1}(\phi)_{c}\rightarrow {\cal
% P}_{0}$ 
the map ${\cal T}:{\cal
H}^{-1}(\phi)\rightarrow\widehat{{\cal
 P}}$ 
is injective.
\end{th}
\pro In case $(G,G)$ is simply connected
the fundamental weights are elements in $X(T)$; in particular condition
(\ref{eq:sc}) in lemma \ref{it:l3} is satisfied for every root $\al\in 
R^{+}$ and
our claim follows from lemma \ref{it:l2}. 
As for the case $G$ satisfies condition $b)$, we see from the Dynkin 
diagram of all simple groups of type different from $B_{l}$, $l\geq 1$ and
$G_{2}$  that 
for every $\al\in R^{+}$ there exists another root $\be$ with $<\be,\al>=1$. On 
the other hand the type $G_{2}$ is simply connected.  
\makebox[2cm][r]{$\Box$}
\begin{th} \label{it:cor2}

Let $a\geq 1$ be the cardinality of the subset
$A\subset R^{+}$ of those roots which 
do not satisfy condition {\em (\ref{eq:sc})} in lemma
{\em \ref{it:l3}}. 
%%%% ultima MODIFICA!!!!!!!!!!!!!!!!!!!!!!!!!!!!!!!!!!!!!!!!!!!!!!!!!!!!!!!
%%%%%%%and $l$ the dimension of 
% $span(A)\subset X(T)_{\ze}\otimes{\bf R}$.
If $d$ denotes the degree of $\pi^{*}K$,
the fibre of $\cal T$ consists of at most $2^{a(d-1)}$
points. \end{th}
\pro  
Let $(P,s)\in {\cal H}^{-1}(\phi)$, $\tau(P,s)$ be as in
theorem \ref{it:thw} and suppose there exists a
 pair $(Q,v)\in {\cal H}^{-1}(\phi)$ such that 
$\tau (Q,v)\cong\tau (P,s)$.
Let $\{ t_{hl}\}_{h,l\in  H}$ and $\{ \td{t}_{hl}\}_{h,l\in  H}$  be cocycles 
representing $\tau(P,s)$ and $\tau(Q,v)$ respectively and write
\(
\td{t}_{hl}=\mih^{-1}t_{hl}m_{l} %\label{eq:t1p} 
%\mtj^{-1}\gtji\mti=\ej^{-1}\mij^{-1}\tji\mii \nonumber % \label{eq:gpp}
\)
for suitable holomorphic maps
 $m_{h}:U_{h}\rightarrow T$ with $h\in H$.
From the proof of 
lemma \ref{it:l3} we can assume that for $a$ chosen ramification points 
$q\in {\cal 
D}_{\be}$, one for each $\be\in A$, 
and every other ramification point
$q\in {\cal 
D}_{\be}$ with $\be\notin A$,   
condition $\be(m_{j}(q))=1$ (for suitable $j\in J$)
 holds.
If $(Q,v)$ is distinct from $(P,s)$, by 
  lemmas \ref{it:l2} and
\ref{it:l3} there exists some $\al\in A$
and some $p_{\al}\in U_{j}\cap {\cal D}_{\al}$ (with suitable $j\in J$) 
such that condition
\begin{equation} \label{eq:-1}
\al(\mij(p_{\al}))=-1 \end{equation}
is satisfied. Moreover, two pairs for which relation (\ref{eq:-1}) holds for
exactly the same set of ramification points coincide by  \rim \ref{sec-res}.1.
% Finally, note that for a fixed ramification point $q\in U_{j}\cap {\cal
% D}_{\al}$ with $j\in J$ and $\al\in R^{+}$ we may assume that condition
% $\al(\mij(p))=1$ is always satisfied (since we can multiply each function
% $m_{k}:U_{k}\rightarrow T$ with 
%$k\in H$ by some constant element $e\in T$ with
% $\al(e)=-1$). Thus our claim follows. 
\makebox[2cm][r]{$\Box$} \\[.1cm]
From theorems \ref{it:cor1} and \ref{it:cor2} and from
 proposition \ref{it:prd} we
obtain the following
\begin{cor}
The image under ${\cal F}$ of the generic Hitchin fibre ${\cal H}^{-1}(\phi)$
contains a Zariski open set in ${\cal P}_{0}$.
\end{cor}
 \subsection{The $PGl(2)$ case.}
Let $\phi\in H^{0}(C,K^{2})$ be generic. Let $P$ be a $PGl(2)$-bundle over $C$
and $s\in H^{0}(C,adP\otimes K)$ such that ${\cal H}(P,s)=\phi$. We indicate by
$pr:Gl(2)\rightarrow PGl(2)=Gl(2)/\ci^{*}$ the factor map and as maximal torus
$T\subset PGl(2)$ we choose the one obtained by restricting $pr$ to the maximal
torus $\widetilde{T}\subset Gl(2)$ given by all diagonal matrices. We also set
$\tl=Lie\; T,\; \td{\tl}=Lie\; \widetilde{T}$. In this setting, $\cit
=\phi^{*}(\tl\otimes K)$ is
a ramified double covering of $C$ whose ramification divisor $\cal D$ satisfies
by definition ${\cal O}({\cal D})\cong \pi^{*}K$.

Let $\{ V_{h}\}_{h\in H}$
and $\{ U_{h}\}_{h\in H}$ be open coverings of $C$ and $\cit$ defined as before.
If $\{ g_{hl}:V_{h}\cap V_{l}\rightarrow \pg \} _{h,l\in H}$,  are transition
functions for $P$, it
is known that there exists some rank 2 vector bundle F, hence some principal
$\gl$-bundle $\widetilde{P}$, with transition functions $\td{g}_{hl}$
satisfying \begin{equation} \label{eq:pr}
pr\circ \td{g}_{hl}=g_{hl}\; \; \forall h,l\in H .\end{equation}
Moreover, any
other rank
2 vector bundle $F'$ has the same property if and only if $F'\cong F\otimes L$
for some line bundle $L\in Pic(C)$. Note also that this implies $deg\;
F\equiv deg\; F'\; \; (mod\; \; 2)$ (since $deg(F\otimes L)=deg\; F\cdot deg\;
L^{2}$).
For the sake of simplicity for any $F$ satisfying relation (\ref{eq:pr}) we
write $P=pr(F)$.
For $\widetilde{P}$ as above, we clearly have an isomorphism $ad\;
\widetilde{P}\otimes K\cong (ad\; P\otimes
K)\oplus K$ and given some fixed generic section $x:C\rightarrow K$ we may
define $\td{s}\in H^{0}(ad\; \widetilde{P}\otimes K)$ by $\td{s}=s\oplus x$.
We set $\td{\phi}={\cal H}_{\gl}(\widetilde{P},\td{s})\in H^{0}(C,K\oplus
K^{2})$ ( the subscript indicating that we are in the $\gl$ setting ) and
observe that
the covering $\td{\phi}^{*}(\td{\tl}\otimes K)$ of $C$ coincides with $\cit$.
%
% Let us denote by ${\cal M}(2,n)$ the moduli space of stable rank-2 vector
% bundles of degree $n$. Since the Hitchin map ${\cal H}_{n}:{\cal
% M}(2,n)\rightarrow
% {\cal K}$ is surjective, for each $n\in\ze$ there exists some rank-2 vector
% bundle $F\in {\cal H}_{n}^{-1}(\td{\phi})$ with  $deg\; F=n$ and $h(F)\in {\cal
% H}_{\pg}^{-1}(\phi)$. 
Then it is clear from the argument above 
that we have a surjective map \[
% "h":{\cal H}_{\gl}^{-1}(\td{\phi})_{c}\times \ze\rightarrow {\cal
% H}_{\pg}^{-1}(\phi)\; .\]
% where ${\cal H}_{\gl}^{-1}(\td{\phi})_{c}$ denotes the connected component 
% containing  $(\widetilde{P},\td{s})$ and
%  $\ze$ parametrizes degrees of the  vector bundles associated to the
% natural representation of $\gl$.
"pr":{\cal H}_{\gl}^{-1}(\td{\phi})\rightarrow {\cal
 H}_{\pg}^{-1}(\phi)\; .\]
This also shows that ${\cal H}_{\pg}^{-1}(\phi)$ has two
components ${\cal H}_{\pg}^{-1}(\phi)_{0}$, ${\cal
H}_{\pg}^{-1}(\phi)_{1}$: namely  $(Q,v)\in {\cal H}_{\pg}^{-1}(\phi)$ is
contained in ${\cal H}_{\pg}^{-1}(\phi)_{0}$ or ${\cal
H}_{\pg}^{-1}(\phi)_{1}$ depending on the parity of the degree of those $F$
which satisfy $pr(F)=Q$.

We now look at 
% the "generalized Prym variety" ${\cal P}_{\gl}$ associated to the
our construction in the 
$\gl$ case. If we indicate by $\chi_{1}$ and $\chi_{2}$ the coordinate functions
on $\widetilde{T}$ and set $\td{\al}=\chi_{1}\cdot \chi_{2}^{-1}$,
$\sigma=s_{\td{\al}}$, we have
by definition \[ {\cal P}_{\gl}=\{ Q\otimes \chi '_{1}\oplus \sigma
^{*}Q\otimes \chi '_{2}\mid Q\in J(\cit)\} \equiv J(\cit )\] 
(the one parameter
subgroups $\chi_{i}'$ being defined %as in (\ref{eq:'}) 
by $\chi_{i}(\chi '_{j})=(\chi_{i},\chi_{j}),\; j=1,2$) and
 \[ \widehat{{\cal P}}_{\gl}=Pic(\cit ).\]  
The map 
% ${\cal
% F}_{c}:{\cal
% H}_{\gl}^{-1}(\td{\phi})_{c}\rightarrow J(\cit )$
${\cal
T}:{\cal
H}_{\gl}^{-1}(\td{\phi})\rightarrow Pic(\cit )$ is injective 
(see theorem
\ref{it:cor1}), dominant and 
%is onto, since for each $Q\in Pic(\cit )$ 
%one can recover the
%rank-2 vector bundle $F$ such that $F\times\gl
%pair $(P,s)$ satisfying 
% ${\cal F}_{c}(P,s)=Q\otimes \chi '_{1}\oplus \sigma
% ^{*}Q\otimes \chi '_{2}$ 
%${\cal T}(P,s)=Q\otimes \chi '_{1}\oplus \sigma
%^{*}Q\otimes \chi '_{2}$ 
%from the push down $\pi_{*}Q$ (see \cite{hi} ).
by Hitchin's theory (see \cite{hi} ) it
preserves the parity of the degrees. By the argument above 
the generic fibre of
the map $"pr"$ is a principal homogeneous space
with respect to
\( \Lambda=\{ M\in Pic(\cit)\; \mid M=\pi^{*}L,\; \; L\in Pic(C)\}. \) 
In this setting the map $\pi^{*}:Pic(C)\rightarrow
Pic(\cit)$ is injective
(since $\cit\rightarrow C$ is a ramified covering: see
e.g
\cite{mum} ),
hence
 $\Lambda$
coincides with $Pic(C)$. Since  
$Pic(\cit)^{\mbox{\em \scriptsize even}}/Pic(C)$ and
$Pic(\cit)^{\mbox{\em \scriptsize odd}}/Pic(C)$ are 
both principal homogeneous spaces 
with respect to the connected group $J(\cit)/J(C)$, it follows that the 
components ${\cal H}_{\pg}^{-1}(\phi)_{0}$, ${\cal
H}_{\pg}^{-1}(\phi)_{1}$ are connected.
Now, let  $\chi '$ be the one parameter subgroup in $T\subset\pg$ given
by composing $pr$ with $\chi '_{1}$ (we have $X(T)^{*}=\ze \chi '$).
%, $\al$  being defined by <\al ',\al>=2$).
By definition, we have $\widehat{{\cal P}} _{\pg}={\cal
P}_{\pg}=\{ 
Q\otimes \chi '\; \mid Q\in J(\cit),\sigma^{*}Q\cong Q^{-1}\} $ and,
since $\pi^{*}:J(C)\rightarrow J(\cit)$ is injective,
this is just the Prym variety $P(\cit ,\sigma)\subset J(\cit)$.    
From theorem \ref{it:thw}  the $\widetilde{T}$-bundle 
$\td{\tau}=\tau(\widetilde{P},\td{s})$ has transition functions
$t_{hl}:U_{h}\cap U_{l}\rightarrow 
\widetilde{T}$ of the form
\[ t_{hl}\deta =diag(\; q_{hl}\deta ,\sigma^{*}q_{hl}\deta \cdot
k_{hl}(\pi\deta)\; \;
).\] One
can easily check that the maps \[ pr\circ t_{hl}\deta = q_{hl}\deta
\cdot \sigma^{*}q_{hl}\deta ^{-1}\cdot k_{hl}(\pi\deta)^{-1}:U_{h}\cap
U_{l}\rightarrow \ci^{*}\]
are transition functions for $\tau=\tau(P,s)$. In other words, if we use 
 the
additive notation, we have ${\cal T}_{\pg}(P,s)=(1-\sigma^{*})\circ {\cal
T}_{\gl}(\widetilde{P},\td{s})$.
Moreover, if $\widetilde{P}'$ is another $\gl$-bundle inducing via the factor
map $pr$ the same $\pg$-bundle $P$, we have that $\tau(\widetilde{P}',\td{s})$
has transition functions $t_{hr}\deta\cdot l_{hr}(\pi\deta)$, where $\{
l_{hr}:V_{h}\cap V_{r}\rightarrow \ci^{*}\} _{h,r\in H}$ define some line bundle
$L$ over $C$. We thus have the following commutative diagram:
\[ \begin{array}{rcccl}
%  Pic(\cit) & \! \equiv J(\cit)\times \ze &
% \stackrel{(1-\sigma^{*})}{\longrightarrow }  & P(\cit ,\sigma )  &  \\
% &
% \mbox{\scriptsize $\{ {\cal F}_{n}\}_{n\in\ze}$}\;
% \updownarrow &
% & \uparrow \! \!
% &  \mbox{\scriptsize ${\cal F}_{0}\coprod {\cal F}_{1}$}
% \\
% &  {\cal H}_{\gl }^{-1}(\td{\phi})_{c}\times \ze & 
% \stackrel{"h"}{\longrightarrow }
% & {\cal
% H}_{\pg}^{-1}(\phi)_{0}\coprod
% {\cal H}_{\pg}^{-1}(\phi)_{1}  &  \end{array} \]
 &  Pic(\cit) & 
\stackrel{(1-\sigma^{*})}{\longrightarrow }  & P(\cit ,\sigma )  &  \\
 \mbox{\scriptsize ${\cal T}_{\gl}$} &
\uparrow &
& \uparrow \! \!
&  \mbox{\scriptsize ${\cal T}_{\pg}$}
\\
&  {\cal H}_{\gl }^{-1}(\td{\phi}) & 
\stackrel{"pr"}{\longrightarrow }
& {\cal
H}_{\pg}^{-1}(\phi)_{0}\coprod
{\cal H}_{\pg}^{-1}(\phi)_{1}  &  \end{array} \]
 If we set
%$(1-\sigma^{*}):Pic(\cit)\rightarrow P(\cit ,\sigma)$ (which is surjective by
%definition) is given by
$\Lambda '=\{ N\in Pic(\cit)\; \mid N=\sigma^{*}N\} $,  we see that all
sufficiently general fibres of the dominant map ${\cal T}_{\pg}$
are  principal homogeneous spaces  with respect to 
% ${\cal F}_{0}\coprod {\cal F}_{1}
%${\cal T}:{\cal H}_{\pg}^{-1}(\phi)\rightarrow P(\cit
%,\sigma)$ is surjective and that each of its fibres
%contains as many points as
%are the elements in  
$\Lambda '/\Lambda$. It is known (see
\cite{mum} ) that $\Lambda '/\Lambda $ 
is isomorphic to $(\ze /2\ze )^{(d-1)}$, $d$
being the number
of ramification points of $\cit$ or, in this setting, the degree of $\pi^{*}K$ .
Note here that the number of $\ze /2\ze$ factors reaches its maximum with
respect to the estimate given in theorem \ref{it:cor2} . Since each 
component ${\cal H}_{\pg}^{-1}(\phi)_{c}$ , $c=0,1$, is  connected,
% of the generic Hitchin fibre
  we have that the generic
fibre of ${\cal F}_{c}:{\cal H}_{\pg}^{-1}(\phi)_{c}
\rightarrow P(\cit ,\sigma)$ consists of $2^{(d-2)}$ points.\\[.5cm]
{\em Acknowledgments.}\\
I wish to express my big debt to my advisor Corrado De Concini for sharing
his ideas on the subject and my gratefulness to Vassil Kanev for
discussions of crucial importance concerning the algebro-geometric aspects of
the problem.\\
Also, it is a pleasure for me to thank 
one of the referees 
for his interesting remarks and his contribution in improving the paper.

\end{document}